\documentclass[11pt]{article}

\usepackage{amsmath}
\usepackage{amssymb}
\usepackage{amsthm}
\usepackage{latexsym}
\usepackage{color}
\usepackage{graphicx}
\usepackage{color}

\DeclareSymbolFont{calletters}{OMS}{cmsy}{m}{n}
\DeclareSymbolFontAlphabet{\mathcal}{calletters}

\def\be{\begin{eqnarray}}
\def\ee{\end{eqnarray}}

\def\b*{\begin{eqnarray*}}
\def\e*{\end{eqnarray*}}

\newtheorem{Theorem}{Theorem}[part]
\newtheorem{Definition}{Definition}[part]
\newtheorem{Proposition}{Proposition}[part]

\newtheorem{Assumption}{Assumption}[part]
\newtheorem{Lemma}{Lemma}[part]
\newtheorem{Corollary}{Corollary}[part]
\newtheorem{Remark}{Remark}[part]

\makeatletter \@addtoreset{equation}{section}

\@addtoreset{Definition}{section}

\@addtoreset{Theorem}{section}

\@addtoreset{Proposition}{section}

\@addtoreset{Property}{section}

\@addtoreset{Assumption}{section}

\@addtoreset{Corollary}{section}

\@addtoreset{Lemma}{section}

\@addtoreset{Remark}{section}

\@addtoreset{Example}{section}

\newcommand{\No}[1]{\left\|#1\right\|}     
\newcommand{\abs}[1]{\left|#1\right|}     

\def \H{\mathbb{H}}

\def \P{\mathbb{P}}

\def\Fh{{\hat F}}

\def\Tr#1{{\rm Tr}\left[#1\right]}

\def \Sup{\displaystyle\sup}

\def\einf{{\rm ess \, inf}}
\def\esup{{\rm ess \, sup}}
\def\trace{{\rm Tr}}

\def\={\;=\;}
\def\.{\;.}

\def\eps{\varepsilon}

\def\reff#1{{\rm(\ref{#1})}}

\def\1{{\bf 1}}
\def \ep{\hbox{ }\hfill{ ${\cal t}$~\hspace{-5.1mm}~${\cal u}$   } }

\def \proof{{\noindent \bf Proof. }}
\def \ep{\hbox{ }\hfill$\Box$}

 \def\normeL2#1{\left\|{#1}\right\|_{L^2}}

 \title{Second Order Reflected Backward Stochastic Differential Equations\footnote{Research supported by the Chair {\it Financial Risks} of the {\it Risk Foundation} sponsored by Soci\'et\'e G\'en\'erale, the Chair {\it Derivatives of the Future} sponsored by the {F\'ed\'eration Bancaire Fran\c{c}aise}, and the Chair {\it Finance and Sustainable Development} sponsored by EDF and Calyon.}
}
\author{ Anis {\sc Matoussi} \footnote{University of Maine, Le Mans and CMAP, Ecole Polytechnique, Paris, anis.matoussi@univ-lemans.fr}\and Dylan {\sc Possamai}\footnote{CMAP, Ecole Polytechnique, Paris, dylan.possamai@polytechnique.edu.}      \and Chao {\sc Zhou}\footnote{CMAP, Ecole Polytechnique, Paris, chao.zhou@polytechnique.edu.} }
 \date{\today}

\begin{document}

 \maketitle

\vspace{3mm}

 \begin{abstract} In this article, we build upon the work of Soner, Touzi and Zhang \cite{stz} to define a notion of a second order backward stochastic differential equation reflected on a lower c\`adl\`ag obstacle. We prove existence and uniqueness of the solution under a Lipschitz type assumption on the generator, and we investigate some links between our reflected 2BSDEs and non-classical optimal stopping problems. Finally, we show that reflected 2BSDEs provide a super-hedging price for American options in a market with volatility uncertainty.

\vspace{10mm}

\vspace{5mm}

\end{abstract}
\newpage

\section{Introduction}

Backward stochastic differential equations (BSDEs for short) appeared in Bismut \cite{bis} in the linear case, and then have been widely studied since the seminal paper of Pardoux and Peng \cite{pardpeng}. Their range of applications includes notably probabilistic numerical methods for partial differential equations, stochastic control, stochastic differential games, theoretical economics and financial mathematics. On a filtered probability space $(\Omega,\mathcal F,\left\{\mathcal F_t\right\}_{0\leq t\leq T},\mathbb P)$ generated by an $\mathbb R^d$-valued Brownian motion $B$, a solution to a BSDE consists on finding a pair of progressively measurable processes $(Y,Z)$ such that
$$Y_t=\xi +\int_t^Tf_s(Y_s,Z_s)ds-\int_t^TZ_sdB_s,\text{ } t\in [0,T], \text{ }\mathbb P-a.s.$$
where $f$ (also called the driver) is a progressively measurable function and $\xi$ is an $\mathcal F_T$-measurable random variable.

\vspace{0.5em}
Pardoux and Peng proved existence and uniqueness of the above BSDE provided that the function $f$ is uniformly Lipschitz in $y$ and $z$ and that $\xi$ and $f_s(0,0)$ are square integrable. Reflected backward stochastic differential equations (RBSDEs for short) were introduced by El Karoui, Kapoudjian, Pardoux, Peng and Quenez in \cite{elkarkap}, followed among others by El Karoui, Pardoux and Quenez in \cite{elkarquen} and Bally, Caballero, El Karoui and Fernandez in \cite{bally} to study related obstacle problems for {PDE}'s and American options pricing. In this case, the solution $Y$ of the BSDE is constrained to stay above a given obstacle process $S$. In order to achieve this, a non-decreasing process $K$ is added to the solution
\begin{align*}
&Y_t=\xi +\int_t^Tf_s(Y_s,Z_s)ds-\int_t^TZ_sdB_s+K_T-K_t,\text{ } t\in [0,T], \text{ }\mathbb P-a.s.\\
&Y_t\geq S_t,\text{ } t\in [0,T], \text{ }\mathbb P-a.s.\\
&\int_0^T(Y_s-S_s)dK_s=0, \text{ }\mathbb P-a.s.,
\end{align*}
where the last condition, also known as the Skorohod condition means that the process $K$ is minimal in the sense that it only acts when $Y$ reaches the obstacle $S$. This condition is crucial to obtain the uniqueness of the classical RBSDEs.

\vspace{0.8em}
Following those pioneering works, many authors have tried to relax the assumptions on the driver of the RBSDE and the corresponding obstacle. Hence, Matoussi \cite{mat} and Lepeltier, Matoussi and Xu \cite{lmx} have extended the existence and uniqueness results to generator with arbitrary growth in $y$. Similarly, Hamad\`ene \cite{hama} and Lepeltier and Xu \cite{lx} proved existence and uniqueness when the obstacle is no longer continuous.

\vspace{0.5em}
More recently, motivated by applications in financial mathematics and probabilistic numerical methods for PDEs (see \cite{ftw}), Cheredito, Soner, Touzi and Victoir \cite{cstv} introduced the notion of Second order BSDEs (2BSDEs), which are connected to the larger class of fully nonlinear PDEs. Then, Soner, Touzi and Zhang \cite{stz} provided a complete theory of existence and uniqueness for 2BSDEs under uniform Lipschitz conditions similar to those of Pardoux and Peng. Their key idea was to reinforce the condition that the 2BSDE must hold $\mathbb P-a.s.$ for every probability measure $\mathbb P$ in a non-dominated class of mutually singular measures (see Section \ref{section.1} for precise definitions). In these regards, this theory shares many similarities with the quasi-sure stochastic analysis of Denis and Martini \cite{denis} and the $G$-expectation theory of Peng \cite{peng}.

\vspace{0.5em}
Our aim in this paper is to provide a complete theory of existence and uniqueness of Second order RBSDEs (2RBSDEs) under the Lipschitz-type hypotheses of \cite{stz} on the driver. We will show that in this context, the definition of a 2RBSDE with a lower obstacle $S$ is very similar to that of a 2BSDE. We do not need to add another non-decreasing process, unlike in the classical case, and we do not need to impose a condition similar to the Skorohod condition. The only change necessary is in the minimal condition that the increasing process $K$ of the 2RBSDE must satisfy.

\vspace{0.5em}
The rest of this paper is organized as follows. In Section \ref{section.1}, we recall briefly some notations, provide the precise definition of 2RBSDEs and show how they are connected to classical RBSDEs. Then, in Section \ref{section.2}, we show a representation formula for the solution of a 2RBSDEs which in turn implies uniqueness. We then provide some links between 2RBSDEs and optimal stopping problems. In Section \ref{section.3}, we give a proof of existence by means of r.c.p.d. techniques, as in \cite{posz} for quadratic 2BDSEs. Let us mention that this proof requires to extend existing results on the theory of $g$-martingales of Peng (see \cite{pengg}) to the reflected case. Since to the best of our knowledge, those results do not exist in the literature, we prove them in the Appendix in Section \ref{section.6}. Finally, we use these new objects in Section \ref{section.5} to study the pricing problem of American options in a market with volatility uncertainty.

\section{Preliminaries} \label{section.1}

Let $\Omega:=\left\{\omega\in C([0,T],\mathbb R^d):\omega_0=0\right\}$ be the canonical space equipped with the uniform norm $\No{\omega}_{\infty}:=\sup_{0\leq t\leq T}|\omega_t|$, $B$ the canonical process, $\mathbb P_0$ the Wiener measure, $\mathbb F:=\left\{\mathcal F_t\right\}_{0\leq t\leq T}$ the filtration generated by $B$, and $\mathbb F^+:=\left\{\mathcal F_t^+\right\}_{0\leq t\leq T}$ the right limit of $\mathbb F$. We first recall the notations introduced in \cite{stz}.

\subsection{The Local Martingale Measures}

We will say that a probability measure $\mathbb P$ is a local martingale measure if the canonical process $B$ is a local martingale under $\mathbb P$. By Karandikar \cite{kar}, we know that we can give pathwise definitions of the quadratic variation $\left<B\right>_t$ and its density $\widehat a_t$.

\vspace{0.5em}
Let $\overline{\mathcal P}_W$ denote the set of all local martingale measures $\mathbb P$ such that
\begin{equation}
\left<B\right>_t \text{ is absolutely continuous in $t$ and $\widehat a$ takes values in $\mathbb S_d^{>0}$, }\mathbb P-a.s.
\label{eq:}
\end{equation}
where $\mathbb S_d^{>0}$ denotes the space of all $d\times d$ real valued positive definite matrices.

\vspace{0.5em}
As usual in the theory of $2$BSDEs, we will concentrate on the subclass $\overline{\mathcal P}_s\subset\overline{\mathcal P}_W$ consisting of all probability measures

\begin{equation}
\mathbb P^\alpha:=\mathbb P_0\circ (X^\alpha)^{-1} \text{ where }X_t^\alpha:=\int_0^t\alpha_s^{1/2}dB_s,\text{ }t\in [0,T],\text{ }\mathbb P_0-a.s.
\end{equation}
for some $\mathbb F$-progressively measurable process $\alpha$ taking values in $\mathbb S_d^{>0}$ with $\int_0^T|\alpha_t|dt<+\infty$, $\mathbb {P}_0-a.s.$

\subsection{The non-linear Generator}
We consider a map $H_t(\omega,y,z,\gamma):[0,T]\times\Omega\times\mathbb{R}\times\mathbb{R}^d\times D_H\rightarrow \mathbb{R}$, where $D_H \subset \mathbb{R}^{d\times d}$ is a given subset containing $0$.

\vspace{0.5em}
Define the corresponding conjugate of $H$ w.r.t.$\gamma$ by
\begin{align*}
&F_t(\omega,y,z,a):=\underset{\gamma \in D_H}{\Sup}\left\{\frac12\trace(a\gamma)-H_t(\omega,y,z,\gamma)\right\} \text{ for } a \in \mathbb S_d^{>0},\\[0.3em]
&\widehat{F}_t(y,z):=F_t(y,z,\widehat{a}_t) \text{ and } \widehat{F}_t^0:=\widehat{F}_t(0,0).
\end{align*}

\vspace{0.5em}
We denote by $D_{F_t(y,z)}:=\left\{a, \text{ }F_t(\omega,y,z,a)<+\infty\right\}$ the domain of $F$ in $a$ for a fixed $(t,\omega,y,z)$.

\vspace{0.5em}
As in \cite{stz} we fix a constant $\kappa \in (1,2]$ and restrict the probability measures in $\mathcal{P}_H^\kappa\subset \overline{\mathcal{P}}_S$

\begin{Definition}\label{def}
$\mathcal{P}_H^\kappa$ consists of all $\mathbb P \in \overline{\mathcal{P}}_S$ such that
$$\underline{a}^\mathbb P\leq \widehat{a}\leq \bar{a}^\mathbb P, \text{ } dt\times d\mathbb P-a.s. \text{ for some } \underline{a}^\P, \bar{a}^\P \in \mathbb S_d^{>0}, \text{ and } \mathbb{E}^{\mathbb{P}}\left[\left(\int_0^T\abs{\widehat{F}_t^0}^\kappa dt\right)^{\frac2\kappa}\right]<+\infty$$

\end{Definition}

\vspace{0.5em}
\begin{Definition}
We say that a property holds $\mathcal{P}^\kappa_H$-quasi-surely ($\mathcal{P}^\kappa_H$-q.s. for short) if it holds $\mathbb P $-a.s. for all $\mathbb P\in \mathcal{P}^\kappa_H$.
\end{Definition}

We now state our main assumptions on the function $F$ which will be our main interest in the sequel
\begin{Assumption} \label{assump.href}
\begin{itemize}
\item[\rm{(i)}] The domain $D_{F_t(y,z)}=D_{F_t}$ is independent of $(\omega,y,z)$.
\item[\rm{(ii)}] For fixed $(y,z,a)$, $F$ is $\mathbb{F}$-progressively measurable in $D_{F_t}$.
\item[\rm{(iii)}] We have the following uniform Lipschitz-type property in $y$ and $z$
$$\forall (y,y',z,z',t,a,\omega), \text{ } \abs{ F_t(\omega,y,z,a)- F_t(\omega,y',z',a)}\leq C\left(\abs{y-y'}+\abs{ a^{1/2}\left(z-z'\right)}\right).$$
\item[\rm{(iv)}] $F$ is uniformly continuous in $\omega$ for the $||\cdot||_\infty$ norm.
\end{itemize}
\end{Assumption}

\vspace{0.5em}
\begin{Remark}
The assumptions (i) and (ii) are classic in the second order framework (\cite{stz}). The Lipschitz assumption (iii) is standard in the BSDE theory since the paper \cite{pardpeng}. The last hypothesis (iv) is also proper to the second order framework, it is linked to our intensive use of regular conditional probability distributions (r.c.p.d.) in our existence proof, and to the fact that we construct our solutions pathwise, thus avoiding complex issues related to negligible sets.
\end{Remark}

\begin{Remark}
  \begin{enumerate}
   \item[$\rm{(i)}$] $\mathcal{P}^\kappa_H$ is decreasing in $\kappa$ since for $\kappa_1<\kappa_2$ with H\"older's inequality
$$ \mathbb{E}^{\mathbb{P}}\left[\left(\int_0^T\abs{\widehat{F}_t^0}^{\kappa_1} dt\right)^{\frac{2}{\kappa_1}}\right] \leq C\mathbb{E}^{\mathbb{P}}\left[\left(\int_0^T\abs{\widehat{F}_t^0}^{\kappa_2} dt\right)^{\frac{2}{\kappa_2}}\right].$$
   \item[$\rm{(ii)}$] The Assumption \ref{assump.href}, together with the fact that $\widehat{F}_t^0<+\infty$, $\P$-a.s for every $\P \in \mathcal{P}^\kappa_H$, implies that $\widehat{a}_t \in D_{F_t}$, $dt\times \P$-a.s., for all $\P \in \mathcal{P}^\kappa_H$.
  \end{enumerate}
\end{Remark}

\subsection{The Spaces and Norms}

We now recall from \cite{stz} the spaces and norms which will be needed for the formulation of the second order BSDEs. Notice that all subsequent notations extend to the case $\kappa=1$.

\vspace{0.5em}
For $p\geq 1$, $L^{p,\kappa}_H$ denotes the space of all $\mathcal F_T$-measurable scalar r.v. $\xi$ with
$$\No{\xi}_{L^{p,\kappa}_H}^p:=\underset{\mathbb{P} \in \mathcal{P}_H^\kappa}{\sup}\mathbb E^{\mathbb P}\left[|\xi|^p\right]<+\infty.$$

$\mathbb H^{p,\kappa}_H$ denotes the space of all $\mathbb F^+$-progressively measurable $\mathbb R^d$-valued processes $Z$ with
$$\No{Z}_{\mathbb H^{p,\kappa}_H}^p:=\underset{\mathbb{P} \in \mathcal{P}_H^\kappa}{\sup}\mathbb E^{\mathbb P}\left[\left(\int_0^T|\widehat a_t^{1/2}Z_t|^2dt\right)^{\frac p2}\right]<+\infty.$$

$\mathbb D^{p,\kappa}_H$ denotes the space of all $\mathbb F^+$-progressively measurable $\mathbb R$-valued processes $Y$ with
$$\mathcal P^\kappa_H-q.s. \text{ c\`adl\`ag paths, and }\No{Y}_{\mathbb D^{p,\kappa}_H}^p:=\underset{\mathbb{P} \in \mathcal{P}_H^\kappa}{\sup}\mathbb E^{\mathbb P}\left[\underset{0\leq t\leq T}{\sup}|Y_t|^p\right]<+\infty.$$

$\mathbb I^{p,\kappa}_H$ denotes the space of all $\mathbb F^+$-progressively measurable $\mathbb R$-valued processes $K$ null at $0$ with
$$\mathcal P^\kappa_H-q.s. \text{ c\`adl\`ag and non-decreasing paths, and }\No{K}_{\mathbb I^{p,\kappa}_H}^p:=\underset{\mathbb{P} \in \mathcal{P}_H^\kappa}{\sup}\mathbb E^{\mathbb P}\left[\left(K_T\right)^p\right]<+\infty.$$

For each $\xi \in L^{1,\kappa}_H$, $\mathbb P\in \mathcal P^\kappa_H$ and $t \in [0,T]$ denote
$$\mathbb E_t^{H,\mathbb P}[\xi]:=\underset{\mathbb P^{'}\in \mathcal P^\kappa_H(t^{+},\mathbb P)}{\esup^{\mathbb P}}\mathbb E^{\mathbb P^{'}}_t[\xi] \text{ where } \mathcal P^\kappa_H(t^{+},\mathbb P):=\left\{\mathbb P^{'}\in\mathcal P^\kappa_H:\mathbb P^{'}=\mathbb P \text{ on }\mathcal F_t^+\right\}.$$

Here $\mathbb{E}_t^{\mathbb P}[\xi]:=E^{\mathbb P}[\xi|\mathcal F_t]$. Then we define for each $p\geq \kappa$,
$$\mathbb L_H^{p,\kappa}:=\left\{\xi\in L^{p,\kappa}_H:\No{\xi}_{\mathbb L_H^{p,\kappa}}<+\infty\right\} \text{ where } \No{\xi}_{\mathbb L_H^{p,\kappa}}^p:=\underset{\mathbb P\in\mathcal P^\kappa_H}{\sup}\mathbb E^{\mathbb P}\left[\underset{0\leq t\leq T}{\esup}^{\mathbb P}\left(\mathbb E_t^{H,\mathbb P}[|\xi|^\kappa]\right)^{\frac{p}{\kappa}}\right].$$

Finally, we denote by $\mbox{UC}_b(\Omega)$ the collection of all bounded and uniformly continuous maps $\xi:\Omega\rightarrow \mathbb R$ with respect to the $\No{\cdot}_{\infty}$-norm, and we let $\mathcal L^{p,\kappa}_H$ be the closure of $\mbox{UC}_b(\Omega)$ under the norm $\No{\cdot}_{\mathbb L^{p,\kappa}_H}$, for every $1\leq\kappa \leq p$.

\subsection{Formulation}
First, we consider a process $S$ which will play the role of our lower obstacle. We will always assume that $S$ verifies the following properties
\begin{itemize}
\item[$\rm{(i)}$] $S$ is $\mathbb F$-progressively measurable and c\`adl\`ag.
\item[$\rm{(ii)}$] $S$ is uniformly continuous in $\omega$ in the sense that for all $t$
$$\abs{S_t(\omega)-S_t(\widetilde{\omega})}\leq\rho\left(\No{\omega-\widetilde{\omega}}_t\right),\text{ }\forall \text{ }(\omega,\widetilde{\omega})\in\Omega^2,$$
for some modulus of continuity $\rho$ and where we define $\No{\omega}_t:=\underset{0\leq s\leq t}{\sup}\abs{\omega(s)}$.
\end{itemize}

\vspace{0.5em}
Then, we shall consider the following second order RBSDE ($2$RBSDE for short) with lower obstacle $S$
\begin{equation}
Y_t=\xi +\int_t^T\widehat{F}_s(Y_s,Z_s)ds -\int_t^TZ_sdB_s + K_T-K_t, \text{ } 0\leq t\leq T, \text{  } \mathcal P_H^\kappa-q.s.
\label{2bsderef}
\end{equation}

\vspace{0.5em}
We follow Soner, Touzi and Zhang \cite{stz}. For any $\mathbb{P}\in\mathcal{P}^\kappa_H$, $\mathbb{F}$-stopping time $\tau$, and $\mathcal{F}_\tau$-measurable random variable $\xi \in \mathbb{L}^2(\mathbb P)$, let $(y^\mathbb{P},z^\mathbb{P},k^\mathbb{P}):=(y^\mathbb{P}(\tau,\xi),z^\mathbb{P}(\tau,\xi),k^\mathbb{P}(\tau,\xi))$ denote the unique solution to the following standard RBSDE with obstacle $S$ (existence and uniqueness have been proved under our assumptions by Lepeltier and Xu in \cite{lx})

$$\begin{cases}
\label{rbsde}
&y_t^\mathbb P=\xi+\int_t^\tau\widehat F_s(y_s^\mathbb P,z_s^\mathbb P)ds-\int_t^\tau z_s^\mathbb PdB_s+k_{\tau}^\mathbb{P}-k_{t}^\mathbb{P}, \text{ }0\leq t\leq \tau,\text{ } \mathbb P-a.s.\\
\nonumber&y_t^\mathbb P\geq S_t, \text{ }\mathbb P-a.s.\\
\nonumber&\int_0^t\left(y^\mathbb P_{s^-}-S_{s^-}\right)dk^\mathbb P_s=0, \text{ }\mathbb P-a.s., \text{ }\forall t\in[0,T].
\end{cases}$$

\vspace{0.5em}
\begin{Definition}
For $\xi \in \mathbb L^{2,\kappa}_H$, we say $(Y,Z)\in \mathbb D^{2,\kappa}_H\times\mathbb H^{2,\kappa}_H$ is a solution to the 2RBSDE \reff{2bsderef} if

\begin{itemize}
\item[$\bullet$] $Y_T=\xi$, and $Y_t\geq S_t$, $t\in[0,T]$, $\mathcal P_H^\kappa-q.s$.
\item[$\bullet$] $\forall \mathbb P \in \mathcal P_H^\kappa$, the process $K^{\mathbb P}$ defined below has non-decreasing paths $\mathbb P-a.s.$
\begin{equation}
K_t^{\mathbb P}:=Y_0-Y_t - \int_0^t\widehat{F}_s(Y_s,Z_s)ds+\int_0^tZ_sdB_s, \text{ } 0\leq t\leq T, \text{  } \mathbb P-a.s.
\label{2bsde.kref}
\end{equation}

\item[$\bullet$] We have the following minimum condition
\begin{equation}
K_t^\mathbb P-k_t^\mathbb P=\underset{ \mathbb{P}^{'} \in \mathcal{P}_H(t^+,\mathbb{P}) }{ \einf^{\mathbb P} }\mathbb{E}_t^{\mathbb P^{'}}\left[K_T^{\mathbb P^{'}}-k_T^{\mathbb P^{'}}\right], \text{ } 0\leq t\leq T, \text{  } \mathbb P-a.s., \text{ } \forall \mathbb P \in \mathcal P_H^\kappa.
\label{2bsde.minref}
\end{equation}
\end{itemize}
\end{Definition}

\vspace{0.8em}
\begin{Remark}
In our proof of existence, we will actually show, using recent results of Nutz \cite{nutz}, that under additional assumptions (related to axiomatic set theory) the family $\left(K^\mathbb P\right)_{\mathbb P\in\mathcal P^\kappa_H}$ can always be aggregated into a universal process $K$.
\end{Remark}

\vspace{0.5em}
Following \cite{stz}, in addition to Assumption \ref{assump.href}, we will always assume

\begin{Assumption}\label{assump.h2ref}
\begin{itemize}
\item[\rm{(i)}] $\mathcal P_H^\kappa$ is not empty.
\item[\rm{(ii)}] The processes $\widehat F^0 $ and $S$ satisfy the following integrability conditions
\begin{align}
\phi^{2,\kappa}_H&:=\underset{\mathbb P\in\mathcal P^\kappa_H}{\sup}\mathbb E^{\mathbb P}\left[\underset{0\leq t\leq T}{\esup}^{\mathbb P}\left(\mathbb E_t^{H,\mathbb P}\left[\int^T_0|\Fh^0_s|^\kappa ds\right]\right)^{\frac{2}{\kappa}}\right]<+\infty\\
\psi^{2,\kappa}_H&:=\underset{\mathbb P\in\mathcal P^\kappa_H}{\sup}\mathbb E^{\mathbb P}\left[\underset{0\leq t\leq T}{\esup}^{\mathbb P}\left(\mathbb E_t^{H,\mathbb P}\left[\left(\underset{0\leq s\leq T}{\sup}\left(S_s\right)^{+}\right)^{\kappa}\right]\right)^{\frac{2}{\kappa}}\right]<+\infty.
\end{align}
\end{itemize}
\end{Assumption}

\subsection{Connection with standard RBSDEs}

If $H$ is linear in $\gamma$, that is to say
$$H_t(y,z,\gamma):=\frac12\Tr{a_t^0\gamma}-f_t(y,z),$$
where $a^0:[0,T]\times\Omega\rightarrow \mathbb S_d^{>0}$ is $\mathbb F$-progressively measurable and has uniform upper and lower bounds. As in \cite{stz}, we no longer need to assume any uniform continuity in $\omega$ in this case. Besides, the domain of $F$ is restricted to $a^0$ and we have
$$\widehat F_t(y,z)=f_t(y,z).$$

If we further assume that there exists some $\mathbb P\in \overline{\mathcal P}_S$ such that $\widehat a$ and $a^0$ coincide $\mathbb P-a.s.$ and $\mathbb E^\mathbb P\left[\int_0^T\abs{f_t(0,0)}^2dt\right]<+\infty$, then $\mathcal P^\kappa_H=\left\{\mathbb P\right\}$.

\vspace{0.5em}
Then, unlike with $2$BSDEs, it is not immediate from the minimum condition \reff{2bsde.minref} that the process $K^\mathbb P-k^\mathbb P$ is actually null. However, we know that $K^\mathbb P-k^\mathbb P$ is a martingale with finite variation. Since $\mathbb P$ satisfy the martingale representation property, this martingale is also continuous, and therefore it is null. Thus we have
$$0=k^\mathbb P-K^\mathbb P,\text{ }\mathbb P-a.s.,$$
and the $2$RBSDE is equivalent to a standard RBSDE. In particular, we see that the part of $K^\mathbb P$ which increases only when $Y_{t^-}>S_{t^-}$ is null, which means that $K^\mathbb P$ satisfies the usual Skorohod condition with respect to the obstacle.

\section{Uniqueness of the solution and other properties}\label{section.2}

\subsection{Representation and uniqueness of the solution}

We have similarly as in Theorem $4.4$ of \cite{stz}

\begin{Theorem}\label{uniqueref}
Let Assumptions \ref{assump.href} and \ref{assump.h2ref} hold. Assume $\xi \in \mathbb{L}^{2,\kappa}_H$ and that $(Y,Z)$ is a solution to $2$RBSDE \reff{2bsderef}. Then, for any $\mathbb{P}\in\mathcal{P}^\kappa_H$ and $0\leq t_1< t_2\leq T$,
\begin{align}
\label{representationref}
Y_{t_1}&=\underset{\mathbb{P}^{'}\in\mathcal{P}^\kappa_H(t_1^+,\mathbb{P})}{\esup^\mathbb{P}}y_{t_1}^{\mathbb{P}^{'}}(t_2,Y_{t_2}), \text{ }\mathbb{P}-a.s.
\end{align}

Consequently, the $2$RBSDE \reff{2bsderef} has at most one solution in $ \mathbb D^{2,\kappa}_H\times\mathbb H^{2,\kappa}_H$.
\end{Theorem}

\begin{Remark}\label{rem.min}
Let us now justify the minimum condition \reff{2bsde.minref}. Assume for the sake of clarity that the generator $\widehat F$ is equal to $0$. By the above Theorem, we know that if there exists a solution to the $2$RBSDE \reff{2bsderef}, then the process $Y$ has to satisfy the representation \reff{representationref}. Therefore, we have a natural candidate for a possible solution of the $2$RBSDE. Now, assume that we could construct such a process $Y$ satisfying the representation \reff{representationref} and which has the decomposition \reff{2bsderef}. Then, taking conditional expectations in $Y-y^\mathbb P$, we end up with exactly the minimum condition \reff{2bsde.minref}.
\end{Remark}

\proof
The proof follows the lines of the proof of Theorem $4.4$ in \cite{stz}.
\vspace{0.5em}
We first assume that \reff{representationref} is true, then 
$$Y_{t}=\underset{\mathbb{P}^{'}\in\mathcal{P}^\kappa_H(t^+,\mathbb{P})}{\esup^\mathbb{P}}y_{t}^{\mathbb{P}^{'}}(T,\xi), \text{ } t\in [0,T], \text{ }\mathbb{P}-a.s., \text{ for all }\mathbb P\in \mathcal{P}^\kappa_H,$$
and thus $Y$ is unique. Since we have that $d\left<Y,B\right>_t=Z_td\left<B\right>_t, \text{ } \mathcal{P}^\kappa_H-q.s.$, $Z$ is unique. Finally, the process $K^\mathbb P$ is uniquely determined. We shall now prove \reff{representationref}.

\begin{itemize}
\item[\rm{(i)}] Fix $0\leq t_1<t_2\leq T$ and $\mathbb P\in\mathcal P^\kappa_H$. For any $\mathbb P^{'}\in\mathcal P^\kappa_H(t_1^+,\mathbb P)$, we have
\end{itemize}
$$Y_{t} = Y_{t_2} + \int_{t}^{t_2} \widehat{F}_s(Y_s,Z_s)ds - \int_{t}^{t_2} Z_sdB_s + K_{t_2}^{\mathbb P^{'}}-K_{t}^{\mathbb P^{'}}, \text{ } t_1\leq t\leq t_2, \text{ } \mathbb P^{'}-a.s.$$

Now, it is clear that we can always decompose the non-decreasing process $K^\mathbb P$ into
$$K^{\mathbb P^{'}}_t=A_t^{\mathbb P^{'}}+B_t^{\mathbb P^{'}},\text{ }\mathbb P^{'}-a.s.,$$
where $A^{\mathbb P^{'}}$ and $B^{\mathbb P^{'}}$ are two non-decreasing processes such that $A^{\mathbb P^{'}}$ only increases when $Y_{t^-}=S_{t^-}$ and $B^{\mathbb P^{'}}$ only increases when $Y_{t^{-}}>S_{t^-}$. With that decomposition, we can apply a generalization of the usual comparison theorem proved by El Karoui et al. \cite{elkarkap}, whose proof is postponed to the appendix, under $\mathbb P^{'}$ to obtain $Y_{t_1}\geq y_{t_1}^{\mathbb P^{'}}(t_2,Y_{t_2})$ and $A_{t_2}^{\mathbb P^{'}}-A_{t_1}^{\mathbb P^{'}}\leq k^{\P^{'}}_{t_2}-k^{\P^{'}}_{t_1},\ \mathbb P^{'}-a.s.$ Since $\mathbb P^{'}=\mathbb P$ on $\mathcal F_t^+$, we get $Y_{t_1}\geq y_{t_1}^{\mathbb P^{'}}(t_2,Y_{t_2})$, $\mathbb P-a.s.$ and thus
$$Y_{t_1}\geq\underset{\mathbb{P}^{'}\in\mathcal{P}^\kappa_H(t_1^+,\mathbb{P})}{\esup^\mathbb{P}}y_{t_1}^{\mathbb{P}^{'}}(t_2,Y_{t_2}), \text{ }\mathbb{P}-a.s.$$

\begin{itemize}
\item[\rm{(ii)}] We now prove the reverse inequality. Fix $\mathbb P\in\mathcal P^\kappa_H$. We will show in $\rm{(iii)}$ below that
\end{itemize}
$$C_{t_1}^{\mathbb P}:=\underset{\mathbb{P}^{'}\in\mathcal{P}^\kappa_H(t_1^+,\mathbb{P})}{\esup^\mathbb{P}}\mathbb E_{t_1}^{\mathbb{P}^{'}}\left[\left(K_{t_2}^{\mathbb P^{'}}-k_{t_2}^{\mathbb P^{'}}-K_{t_1}^{\mathbb P^{'}}+k_{t_1}^{\mathbb P^{'}}\right)^2\right]<+\infty,\text{ }\mathbb P-a.s.$$

\vspace{0.5em}
For every $\mathbb P^{'}\in \mathcal P^\kappa_H(t^+,\mathbb P)$, denote
$$\delta Y:=Y-y^{\mathbb P^{'}}(t_2,Y_{t_2}),\text{ }\delta Z:=Z-z^{\mathbb P^{'}}(t_2,Y_{t_2}) \text{ and } \delta K^{\mathbb P^{'}}:=K^{\mathbb P^{'}}-k^{\mathbb P^{'}}(t_2,Y_{t_2}).$$

By the Lipschitz Assumption \ref{assump.href}$\rm{(iii)}$ and using a classical linearization procedure, we can define a continuous process $M$ such that for all $p\geq 1$

\begin{equation}
\label{truc.M}
\mathbb E_{t_1}^{\mathbb P^{'}}\left[\underset{t_1\leq t\leq t_2}{\sup}(M_t)^p+\underset{t_1\leq t\leq t_2}{\sup}(M_t^{-1})^p\right]\leq C_p,\text{ }\mathbb P^{'}-a.s.,
\end{equation}
and
\begin{equation}
\label{brubru}
\delta Y_{t_1}=\mathbb E_{t_1}^{\mathbb P^{'}} \left[\int_{t_1}^{t_2}M_{t^-} d\delta K_t^{\mathbb P^{'}}\right].
\end{equation}

\vspace{0.5em}
Let us now prove that the process $K^{\mathbb P^{'}}-k^{\mathbb P^{'}}$ is non-decreasing. By the minimum condition \reff{2bsde.minref}, it is clear that it is actually a $\mathbb P^{'}$-submartingale. Let us apply the Doob-Meyer decomposition under $\mathbb P^{'}$, we get the existence of a $\mathbb P^{'}$-martingale $N^{\mathbb P^{'}}$ and a non-decreasing process $P^{\mathbb P^{'}}$, both null at $0$, such that
$$K_t^{\mathbb P^{'}}-k_t^{\mathbb P^{'}}=N_t^{\mathbb P^{'}}+P_t^{\mathbb P^{'}},\text{ }\mathbb P^{'}-a.s.$$

Then, since we know that all the probability measures in $\mathcal P^\kappa_H$ satisfy the martingale representation property, the martingale $N^{\mathbb P^{'}}$ is continuous. Besides, by the above equation, it also has finite variation. Hence, we have $N^{\mathbb P^{'}}=0$, and the result follows. Returning back to \reff{brubru}, we can now write
\begin{align*}
\delta Y_{t_1}&\leq\mathbb E_{t_1}^{\mathbb P^{'}}\left[\underset{t_1\leq t\leq t_2}{\sup}(M_t)\left(\delta K_{t_2}^{\mathbb P^{'}}-\delta K_{t_1}^{\mathbb P^{'}}\right)\right]\\
&\leq \left(\mathbb E_{t_1}^{\mathbb P^{'}}\left[\underset{t_1\leq t\leq t_2}{\sup}(M_t)^3\right]\right)^{1/3}\left(\mathbb E_{t_1}^{\mathbb P^{'}}\left[\left(\delta K_{t_2}^{\mathbb P^{'}}-\delta K_{t_1}^{\mathbb P^{'}}\right)^{3/2}\right]\right)^{2/3}\\
&\leq \left(\mathbb E_{t_1}^{\mathbb P^{'}}\left[\underset{t_1\leq t\leq t_2}{\sup}(M_t)^3\right]\right)^{1/3}\left(\mathbb E_{t_1}^{\mathbb P^{'}}\left[\delta K_{t_2}^{\mathbb P^{'}}-\delta K_{t_1}^{\mathbb P^{'}}\right]\mathbb E_{t_1}^{\mathbb P^{'}}\left[\left(\delta K_{t_2}^{\mathbb P^{'}}-\delta K_{t_1}^{\mathbb P^{'}}\right)^{2}\right]\right)^{1/3}\\
&\leq C(C_{t_1}^\mathbb P)^{1/3}\left(\mathbb E_{t_1}^{\mathbb P^{'}}\left[\delta K_{t_2}^{\mathbb P^{'}}-\delta K_{t_1}^{\mathbb P^{'}}\right]\right)^{1/3},\text{ }\mathbb P-a.s.
\end{align*}

By taking the essential infimum in $\mathbb P^{'}\in\mathcal P^\kappa_H(t_1^+,\mathbb P)$ on both sides and using the minimum condition \reff{2bsde.minref}, we obtain the reverse inequality.

\begin{itemize}
\item[\rm{(iii)}] It remains to show that the estimate for $C_{t_1}^{\mathbb P}$ holds. But by definition, we clearly have
\end{itemize}
\begin{align*}
\mathbb E^{\mathbb{P}^{'}}\left[\left(K_{t_2}^{\mathbb P^{'}}-k_{t_2}^{\mathbb P^{'}}-K_{t_1}^{\mathbb P^{'}}+k_{t_1}^{\mathbb P^{'}}\right)^2\right]&\leq C\left(\No{Y}^2_{\mathbb D^{2,\kappa}_H}+\No{Z}^{2}_{\H^{2,\kappa}_H}+\phi^{2,\kappa}_H\right)\\
&\hspace{0.9em}+C\underset{\mathbb P\in\mathcal P^\kappa_H}{\sup}\mathbb E^\mathbb P\left[\underset{0\leq t\leq T}{\sup}\abs{y_t^\mathbb P}^2+\int_0^T\abs{\widehat a_t^{1/2}z_s^\mathbb P}^2ds\right]<+\infty,
\end{align*}
since the last term on the right-hand side is finite thanks to the integrability assumed on $\xi$ and $\widehat F^0$. Then we can proceed exactly as in the proof of Theorem $4.4$ in \cite{stz}.
\ep

\vspace{0.5em}
Finally, the following comparison Theorem follows easily from the classical one for RBSDEs (see for instance Theorem $3.4$ in \cite{lx}) and the representation \reff{representationref}.

\begin{Theorem}
Let $(Y,Z)$ and $(Y',Z')$ be the solutions of $2$RBSDEs with terminal conditions $\xi$ and $\xi^{'}$, lower obstacles $S$ and $S^{'}$ and generators $\widehat F$ and $\widehat F^{'}$ respectively (with the corresponding functions $H$ and $H^{'}$), and let $(y^\mathbb P,z^\mathbb P,k^\mathbb P)$ and $(y'^{\mathbb P},z'^{\mathbb P},k'^{\mathbb P})$ the solutions of the associated RBSDEs. Assume that they both verify our Assumptions \ref{assump.href} and \ref{assump.h2ref}, that $\mathcal P^\kappa_H\subset\mathcal P^\kappa_{H^{'}}$ and that we have $\mathcal P^\kappa_H-q.s.$
$$\xi\leq\xi^{'}, \ \widehat F_t(y'^{\mathbb P}_t,z'^{\mathbb P}_t)\leq \widehat F^{'}_t(y'^{\mathbb P}_t,z'^{\mathbb P}_t),\text{ and } S_t\leq S_t^{'}.$$

Then $Y\leq Y'$, $\mathcal P^\kappa_H-q.s.$
\end{Theorem}

\begin{Remark}
Note that in our context, in the above comparison Theorem, even if the obstacles $S$ and $S^{'}$ are identical, we cannot compare the increasing processes $K^\mathbb P$ and $K'^{\mathbb P}$. This is due to the fact that the processes $K^\mathbb P$ do not satisfy the Skorohod condition, since it can be considered, at least formally, to come from the addition of an increasing process due to the fact that we work with second-order BSDEs, and an increasing process due to the reflection constraint. And only the second one is bound to satisfy the Skorohod condition.
\end{Remark}

\subsection{Some properties of the solution}

Now that we have proved the representation \reff{representationref}, we can show, as in the classical framework, that the solution $Y$ of the $2$RBSDE is linked to an optimal stopping problem

\begin{Proposition}\label{prop.trr}
Let $(Y,Z)$ be the solution to the above $2$RBSDE \reff{2bsderef}. Then for each $t\in [0,T]$ and for all $\mathbb P\in\mathcal P^\kappa_H$
\begin{align}\label{opt1}
Y_t&=\underset{\mathbb P^{'}\in\mathcal P^\kappa_H(t^+,\mathbb P)}{\esup^\mathbb P}\underset{\tau\in\mathcal T_{t,T}}{\esup}\text{ }\mathbb E_t^{\mathbb P^{'}}\left[\int_t^\tau\widehat F_s(y_s^{\mathbb P^{'}},z_s^{\mathbb P^{'}})ds +S_\tau1_{\tau<T}+\xi1_{\tau=T}\right],\text{ }\mathbb P-a.s.\\[0.8em]
&=\underset{\tau\in\mathcal T_{t,T}}{\esup}\text{ }\mathbb E_t^\mathbb P\left[\int_t^\tau\widehat F_s(
Y_s,Z_s)ds +A_\tau^\mathbb P-A_t^\mathbb P+S_\tau1_{\tau<T}+\xi1_{\tau=T}\right],\text{ }\mathbb P-a.s.
\end{align}
where $\mathcal T_{t,T}$ is the set of all stopping times valued in $[t,T]$ and where $A_t^\mathbb P:=\int_0^t1_{Y_{s^-}>S_{s^-}}dK_s^\mathbb P$ is the part of $K^\mathbb P$ which only increases when $Y_{s^-}>S_{s^-}$.
\end{Proposition}

\vspace{0.5em}

\begin{Remark}\label{rem.low}
We want to highlight here that unlike with classical RBSDEs, considering a lower obstacle in our context is fundamentally different from considering an upper obstacle. Indeed, having an lower obstacle corresponds, at least formally, to add an increasing process in the definition of a 2BSDE. Since there is already an increasing process in that definition, we still end up with an increasing process. However, in the case of a upper obstacle, we would have to add a decreasing process in the definition, therefore ending up with a finite variation process. This situation thus becomes much more complicated. Furthermore, in that case we conjecture that the above representation of Proposition \ref{prop.trr} would hold with a sup-inf instead of a sup-sup, indicating that this situation should be closer to stochastic games than to stochastic control. We believe that such a generalization would be extremely interesting from the point of view of applications. Indeed, optimal stopping problems (or cooperative controller-and-stopper games) and zero-sum stochastic controller-and-stopper games (or robust optimal stopping problems) with controlled state process have been actively studied in the literature. To name but a few:

\vspace{0.5em}
Karatzas and Sudderth \cite{ks1} solve an optimal stopping problem in which the controller chooses both the drift coefficient and the volatility coefficient of a linear one-dimensional diffusion along a given interval on $\mathbb R$ and selects a stopping rule to maximize her reward. Under mild regularity conditions, by relying on theorems of
optimal stopping for one-dimensional diffusions, they show that this problem admits a simple solution.

\vspace{0.5em}
In a similar setting, Karatzas and Sudderth \cite{ks2} study a zero-sum stochastic game in which a controller selects the coefficients of a linear diffusion along a given interval on $\mathbb R$ to minimize her cost and a stopper chooses a stopping time to maximize his reward. Under appropriate conditions, they prove that this game has a value and describe fairly explicitly a saddle point of optimal strategies.

\vspace{0.5em}
Bayraktar and Huang \cite{bh} consider a zero-sum stochastic differential controller-and-stopper game in which the
state process is a controlled multi-dimensional diffusion. In this game, while the controller selects both the drift and the volatility terms of the state process to maximize her reward, the stopper chooses a stopping time to minimize his cost. Under appropriate conditions, by proving dynamic-programming-type results, they show that the game has a value and the value function is the unique viscosity solution to an obstacle problem for a Hamilton-Jacobi-Bellman equation. Their results can also be interpreted as a solution to a robust optimal stopping problem under both drift and volatility uncertainty.

\vspace{0.5em}
We also refer the reader to Karatzas and Zamfirescu \cite{kz}, Bayraktar, Karatzas and Yao \cite{bky}, Bayraktar and Yao \cite{by1},\cite{by2} among others, for the case where there is only drift uncertainty.

\vspace{0.5em}
We believe that the theory of 2RBSDEs could provide interesting new tools to tackle the above problems or their possible extensions.
\end{Remark}

\vspace{0.5em}
\proof
By Proposition $3.1$ in \cite{lx}, we know that for all $\mathbb P\in\mathcal P^\kappa_H$
$$y_t^\mathbb P=\underset{\tau\in\mathcal T_{t,T}}{\esup}\text{ }\mathbb E_t^{\mathbb P}\left[\int_t^\tau\widehat F_s(y_s^{\mathbb P},z_s^{\mathbb P})ds +S_\tau1_{\tau<T}+\xi1_{\tau=T}\right],\text{ }\mathbb P-a.s.$$

Then the first equality is a simple consequence of the representation formula \reff{representationref}. For the second one, we proceed exactly as in the proof of Proposition $3.1$ in \cite{lx}. Fix some $\mathbb P\in\mathcal P^\kappa_H$ and some $t\in[0,T]$. Let $\tau\in\mathcal T_{t,T}$. We obtain by taking conditional expectation in \reff{2bsderef}
\begin{align*}
Y_t&=\mathbb E_t^{\mathbb P}\left[Y_\tau+\int_t^\tau\widehat F_s(Y_s,Z_s)ds+K^{\mathbb P}_\tau-K^{\mathbb P}_t\right]\\
&\geq\mathbb E_t^{\mathbb P}\left[\int_t^\tau\widehat F_s(Y_s,Z_s)ds+S_\tau1_{\tau<T}+\xi1_{\tau=T}+A^{\mathbb P}_\tau-A^{\mathbb P}_t\right].
\end{align*}

\vspace{0.5em}
This implies that
$$Y_t\geq\underset{\tau\in\mathcal T_{t,T}}{\esup}\text{ }\mathbb E_t^\mathbb P\left[\int_t^\tau\widehat F_s(
Y_s,Z_s)ds +A_\tau^\mathbb P-A_t^\mathbb P+S_\tau1_{\tau<T}+\xi1_{\tau=T}\right],\text{ }\mathbb P-a.s.$$

\vspace{0.5em}
Fix some $\eps>0$ and define the stopping time $D^{\mathbb P,\eps}_t:=\inf\left\{u\geq t,\text{ }Y_u\leq S_u+\eps,\text{ }\mathbb P-a.s.\right\}\wedge T$. It is clear by definition that on the set $\left\{D_t^{\mathbb P,\eps}<T\right\}$, we have $Y_{D^{\mathbb P,\eps}_t}\leq S_{D^{\mathbb P,\eps}_t}+\eps$. Similarly, on the set $\left\{D^{\mathbb P,\eps}_t=T\right\}$, we have $Y_s>S_s+\eps$, for all $t\leq s\leq T$. Hence, for all $s\in[t,D^{\mathbb P,\eps}_t]$, we have $Y_{s^-}>S_{s^-}$. This implies that $K_{D^{\mathbb P,\eps}_t}-K_t=A_{D^{\mathbb P,\eps}_t}-A_t$, and therefore

$$Y_t\leq\mathbb E_t^\mathbb P\left[\int_t^{D^{\mathbb P,\eps}_t}\widehat F_s(
Y_s,Z_s)ds +A_{D^{\mathbb P,\eps}_t}^\mathbb P-A_t^\mathbb P+S_{D^{\mathbb P,\eps}_t}1_{D^{\mathbb P,\eps}_t<T}+\xi1_{D^{\mathbb P,\eps}_t=T}\right]+\eps,$$
which ends the proof by arbitrariness of $\eps$.
\ep


We now show that we can obtain more information about the non-decreasing processes $K^\mathbb P$.

\begin{Proposition}\label{prop.imp}
Let Assumptions \ref{assump.href} and \ref{assump.h2ref} hold. Assume $\xi\in\mathbb L^{2,\kappa}_H$ and $(Y,Z)\in \mathbb D^{2,\kappa}_H\times\mathbb H^{2,\kappa}_H$ is a solution to the 2RBSDE \reff{2bsderef}. Let $\left\{(y^\mathbb P,z^\mathbb P,k^\mathbb P)\right\}_{\mathbb P\in\mathcal P^\kappa_H}$ be the solutions of the corresponding BSDEs \reff{rbsde}. Then we have the following result. For all $t\in[0,T]$,
$$\int_0^t1_{\{Y_{s^-}=S_{s^-}\}}dK_s^\mathbb P=\int_0^t1_{\{Y_{s^-}=S_{s^-}\}}dk_s^\mathbb P, \text{ $\mathbb P-a.s.$}$$
\end{Proposition}

\vspace{0.5em}
\proof
Let us fix a given $\mathbb P\in\mathcal P^\kappa_H$. Let $\tau_1$ and $\tau_2$ be two $\mathbb P$-stopping times such that for all $t\in[\tau_1,\tau_2)$, $Y_{t^-}=S_{t^-}$, $\mathbb P-a.s.$

\vspace{0.5em}
First, by the representation formula \reff{representationref}, we necessarily have for all $\mathbb P$, $Y_{t^-}\geq y^\mathbb P_{t^-}$, $\mathbb P-a.s.$ for all $t$. Moreover, since we also have $y^\mathbb P_t\geq S_t$ by definition, this implies, since all the processes here are c\`adl\`ag, that we must have
$$Y_{t^-}=y^\mathbb P_{t^-}=S_{t^-},\ t\in[\tau_1,\tau_2),\ \mathbb P-a.s.$$

Using the fact that $Y$ and $y^\mathbb P$ solve respectively a 2BSDE and a BSDE, we also have
$$S_{t^{-}}+\Delta Y_t=Y_t=Y_{u}+\int_t^{u}\widehat F_s(Y_s,Z_s)ds-\int_t^{u}Z_sdB_s+K_{u}^\mathbb P-K_t^{\mathbb P},\ \tau_1\leq t\leq u< \tau_2,\ \mathbb P-a.s.,$$
and
$$S_{t^{-}}+\Delta y^\mathbb P_t=Y_t=y^\mathbb P_{u}+\int_t^{u}\widehat F_s(y^\mathbb P_s,z^\mathbb P_s)ds-\int_t^{u}z^\mathbb P_sdB_s+k_{u}^\mathbb P-k_t^{\mathbb P},\ \tau_1\leq t\leq u< \tau_2,\ \mathbb P-a.s.$$

Identifying the martingale parts above, we obtain that $Z_s=z^\mathbb P_s$, $\mathbb P-a.s.$ for all $s\in [t,u]$. Then, identifying the finite variation parts, we have
\begin{align*}
&\Delta Y_u-\Delta Y_t+\int_t^u\widehat F_s(Y_s,Z_s)ds+K_{u}^\mathbb P-K_t^{\mathbb P}=\Delta y^\mathbb P_u-\Delta y^\mathbb P_t+\int_t^u\widehat F_s(y^\mathbb P_s,z^\mathbb P_s)ds+k_{u}^\mathbb P-k_t^{\mathbb P}.
\end{align*}

Now, we clearly have
$$\int_t^u\widehat F_s(Y_s,Z_s)ds=\int_t^u\widehat F_s(y^\mathbb P_s,z^\mathbb P_s)ds,$$
since $Z_s=z^\mathbb P_s$, $\mathbb P-a.s.$  and $Y_{s^-}=y^\mathbb P_{s^-}=S_{s^-}$ for all $s\in [t,u]$. Moreover, since $Y_{s^-}=y^\mathbb P_{s^-}=S_{s^-}$ for all $s\in [t,u]$ and since all the processes are c\`adl\`ag, the jumps of $Y$ and $y^\mathbb P$ are equal to the jumps of $S$. Therefore, we can further identify the finite variation part to obtain
$$K_{u}^\mathbb P-K_t^{\mathbb P}=k_{u}^\mathbb P-k_t^{\mathbb P},$$
which is the desired result.
\ep

\vspace{0.5em}
\begin{Remark}
Recall that at least formally, the role of the non-decreasing processes $K^\mathbb P$ is on the one hand to keep the solution of the 2RBSDE above the obstacle $S$ and on the other hand to keep it above the corresponding RBSDE solutions $y^\mathbb P$, as confirmed by the representation formula \reff{representationref}. What the above result tells us is that if $Y$ becomes equal to the obstacle, then it suffices to push it exactly as in the standard RBSDE case. This is conform to the intuition. Indeed, when $Y$ reaches $S$, then all the $y^\mathbb P$ are also on the obstacle, therefore, there is no need to counter-balance the second order effects.
\end{Remark}

\begin{Remark}\label{rem.decomp}
The above result leads us naturally to think that one could decompose the non-decreasing process $K^\mathbb P$ into two non-decreasing processes $A^\mathbb P$ and $V^\mathbb P$ such that $A^\mathbb P$ satisfies the usual Skorohod condition and $V^\mathbb P$ satisfies
\begin{equation*}
V_t^\mathbb P=\underset{ \mathbb{P}^{'} \in \mathcal{P}_H^\kappa(t^+,\mathbb{P}) }{ \einf^{\mathbb P} }\mathbb{E}_t^{\mathbb P^{'}}\left[V_T^{\mathbb P^{'}}\right], \text{ } 0\leq t\leq T, \text{  } \mathbb P-a.s., \text{ } \forall \mathbb P \in \mathcal P_H^\kappa.
\end{equation*}

Such a decomposition would isolate the effects due to the obstacle and the ones due to the second-order. Of course, the choice $A^\mathbb P:=k^\mathbb P$ would be natural, given the minimum condition \reff{2bsde.minref}. However the situation is not that simple. Indeed, we know that
$$\int_0^t1_{\{Y_{s^-}=S_{s^-}\}}dK_s^\mathbb P=\int_0^t1_{\{Y_{s^-}=S_{s^-}\}}dk_s^\mathbb P.$$

But $k^\mathbb P$ can increase when $Y$ is strictly above the obstacle, since we can have $Y_{t^-}>y_{t^-}^\mathbb P=S_{t^-}$. We can thus only write
$$K^\mathbb P_t=\int_0^t1_{\{Y_{s^-}=S_{s^-}\}}k^\mathbb P_s+V_t^\mathbb P.$$

Then $V^\mathbb P$ satisfies the minimum condition \reff{2bsde.minref} when $Y_{t^-}=S_{t^-}$ and when $y^\mathbb P_{t^-}>S_{t^-}$. However, we cannot say anything when $Y_{t^-}>y_{t^-}^\mathbb P=S_{t^-}$. The existence of such a decomposition, which is also related to the difficult problem of the Doob-Meyer decomposition for the $G$-submartingales of Peng \cite{peng}, is therefore still an open problem.
\end{Remark}

\vspace{0.5em}
As a Corollary of the above result, if we have more information on the obstacle $S$, we can give a more explicit representation for the processes $K^\mathbb P$. The proof comes directly from the above Proposition and Proposition $4.2$ in \cite{elkarquen}.

\begin{Assumption}\label{assump.s}
$S$ is a semi-martingale of the form
$$S_t=S_0+\int_0^tU_sds+\int_0^tV_sdB_s+C_t,\text{ }\mathcal P_H^\kappa-q.s.$$
where $C$ is c\`adl\`ag process of integrable variation such that the measure $dC_t$ is singular with respect to the Lebesgue measure $dt$ and which admits the following decomposition
$$C_t=C_t^+-C_t^-,$$
where $C^+$ and $C^-$ are non-decreasing processes. Besides, $U$ and $V$ are respectively $\mathbb R$ and $\mathbb R^d$-valued $\mathcal F_t$ progressively measurable processes such that
$$\int_0^T(\abs{U_t}+\abs{V_t}^2)dt+C_T^++C_T^- <  +\infty,\text{ }\mathcal P_H^\kappa-q.s.$$
\end{Assumption}

\vspace{0.5em}
\begin{Corollary}
Let Assumptions \ref{assump.href}, \ref{assump.h2ref} and \ref{assump.s} hold. Let $(Y,Z)$ be the solution to the $2$RBSDE \reff{2bsderef}, then
\begin{align}
&Z_t=V_t, \text{ }dt\times\mathcal P^\kappa_H-q.s.\text{ on the set }\left\{Y_{t^-}=S_{t^-}\right\},
\end{align}
and there exists a progressively measurable process $(\alpha_t^\mathbb P)_{0\leq t\leq T}$ such that $0\leq \alpha\leq 1$ and
$$1_{\left\{Y_{t^-}=S_{t^-}\right\}}dK_t^\mathbb P=\alpha_t^\mathbb P1_{\{Y_{s^-}=S_{s^-}\}}\left(\left[\widehat F_t(S_t,V_t)+U_t\right]^-dt+dC_t^-\right).$$
\end{Corollary}

%
%
%
%
%

\subsection{A priori estimates}
We conclude this section by showing some a priori estimates which will prove useful.

\begin{Theorem}\label{estimatesref}
Let Assumptions \ref{assump.href} and \ref{assump.h2ref} hold. Assume $\xi\in\mathbb L^{2,\kappa}_H$ and $(Y,Z)\in \mathbb D^{2,\kappa}_H\times\mathbb H^{2,\kappa}_H$ is a solution to the 2RBSDE \reff{2bsderef}. Let $\left\{(y^\mathbb P,z^\mathbb P,k^\mathbb P)\right\}_{\mathbb P\in\mathcal P^\kappa_H}$ be the solutions of the corresponding RBSDEs \reff{rbsde}. Then, there exists a constant $C_\kappa$ depending only on $\kappa$, $T$ and the Lipschitz constant of $\widehat F$ such that
\begin{align*}
\No{Y}^2_{\mathbb D^{2,\kappa}_H}+\No{Z}^2_{\mathbb H^{2,\kappa}_H}+\underset{\mathbb P\in \mathcal P^\kappa_H}{\sup}\mathbb E^\mathbb P\left[(K_T^\mathbb P)^2\right]&\leq C\left(\No{\xi}^2_{\mathbb L^{2,\kappa}_H}+\phi^{2,\kappa}_H+\psi^{2,\kappa}_H\right)\\
\underset{\mathbb P\in \mathcal P^\kappa_H}{\sup}\left\{\No{y^\mathbb P}^2_{\mathbb D^{2}(\mathbb P)}+\No{z^\mathbb P}^2_{\mathbb H^{2}(\mathbb P)}+\No{k^\mathbb P}^2_{\mathbb I^{2}(\mathbb P)}\right\}&\leq C\left(\No{\xi}^2_{\mathbb L^{2,\kappa}_H}+\phi^{2,\kappa}_H+\psi^{2,\kappa}_H\right).
\end{align*}
\end{Theorem}

\proof
By Lemma $2$ in \cite{ham}, we know that there exists a constant $C_\kappa$ depending only on $\kappa$, $T$ and the Lipschitz constant of $\widehat F$, such that for all $\mathbb P$
\begin{equation}
\label{estimref}
\abs{y_t^\mathbb P}\leq C_\kappa\mathbb E_t^\mathbb P\left[\abs{\xi}^\kappa+\int_t^T\abs{\widehat F^0_s}^\kappa ds+\underset{t\leq s\leq T}{\sup}(S_s^+)^\kappa\right].
\end{equation}

Let us note immediately, that in \cite{ham}, the result is given with an expectation and not a conditional expectation, and more importantly that the process considered are continuous. However, the generalization is easy for the conditional expectation. As far as the jumps are concerned, their proof only uses It\^o's formula for smooth convex functions, for which the jump part can been taken care of easily in the estimates. Then, one can follow exactly their proof to get our result. This immediately provides the estimate for $y^\mathbb P$. Now by definition of our norms, we get from \reff{estimref} and the representation formula \reff{representationref} that
\begin{equation}
\label{eq.y}
\No{Y}_{\mathbb D^{2,\kappa}_H}^2\leq C_\kappa\left(\No{\xi}^2_{\mathbb L^{2,\kappa}_H}+\phi^{2,\kappa}_H+\psi^{2,\kappa}_H\right).
\end{equation}

\vspace{0.5em}
Now apply It\^o's formula to $\abs{Y}^2$ under each $\mathbb P\in\mathcal P^\kappa_H$. We get as usual for every $\epsilon>0$

\begin{align}
\nonumber\mathbb E^\mathbb P\left[\int_0^T\abs{\widehat a_t^{\frac12}Z_t}^2dt\right]&\leq C\mathbb E^\mathbb P\left[\abs{\xi}^2+\int_0^T\abs{Y_t}\left(|\widehat F^0_t|+\abs{Y_t}+|\widehat a_t^{\frac12}Z_t|\right)dt\right]+\mathbb E^\mathbb P\left[\int_0^T\abs{Y_t}dK_t^\mathbb P\right]\\
\nonumber&\leq C\left(\No{\xi}_{\mathbb L^{2,\kappa}_H}+\mathbb E^\mathbb P\left[\underset{0\leq t\leq T}{\sup}\abs{Y_t}^2+\left(\int_0^T\abs{\widehat F^0_t}dt\right)^2\right]\right)\\
&\hspace{0.9em}+\epsilon\mathbb E^\mathbb P\left[\int_0^T\abs{\widehat a_t^{1/2}Z_t}^2dt+\abs{K_T^\mathbb P}^2\right]+\frac{C^2}{\eps}\mathbb E^\mathbb P\left[\underset{0\leq t\leq T}{\sup}\abs{Y_t}^2\right].
\label{grahou}
\end{align}

\vspace{0.5em}
Then by definition of our $2$RBSDE, we easily have
\begin{equation}
\mathbb E^\mathbb P\left[\abs{K_T^\mathbb P}^2\right]\leq C_0\mathbb E^\mathbb P\left[\abs{\xi}^2+\underset{0\leq t\leq T}{\sup}\abs{Y_t}^2+\int_0^T\abs{\widehat a_t^{1/2}Z_t}^2dt+\left(\int_0^T\abs{\widehat F^0_t}dt\right)^2\right],
\label{eq.kk}
\end{equation}
for some constant $C_0$, independent of $\epsilon$.

\vspace{0.5em}
Now set $\epsilon:=(2(1+C_0))^{-1}$ and plug \reff{eq.kk} in \reff{grahou}. One then gets
$$\mathbb E^\mathbb P\left[\int_0^T\abs{\widehat a_t^{1/2}Z_t}^2dt\right]\leq C\mathbb E^\mathbb P\left[\abs{\xi}^2+\underset{0\leq t\leq T}{\sup}\abs{Y_t}^2+\left(\int_0^T\abs{\widehat F^0_t}dt\right)^2\right].$$

\vspace{0.5em}
From this and the estimate for $Y$, we immediately obtain
$$\No{Z}_{\mathbb H^{2,\kappa}_H}\leq C\left(\No{\xi}^2_{\mathbb L^{2,\kappa}_H}+\phi^{2,\kappa}_H+\psi^{2,\kappa}_H\right).$$

\vspace{0.5em}
The estimate for $K^\mathbb P$ comes from \reff{eq.kk} and the ones for $z^\mathbb P$ and $k^\mathbb P$ can be proved similarly.
\ep

\vspace{0.5em}
\begin{Theorem}\label{estimates2}
Let Assumptions \ref{assump.href} and \ref{assump.h2ref} hold. For $i=1,2$, let $(Y^i,Z^i)$ be the solutions to the 2RBSDE \reff{2bsderef} with terminal condition $\xi^i$ and lower obstacle $S$. Then, there exists a constant $C_\kappa$ depending only on $\kappa$, $T$ and the Lipschitz constant of $\widehat F$ such that
\begin{align*}
&\No{Y^1-Y^2}_{\mathbb D^{2,\kappa}_H}\leq C\No{\xi^1-\xi^2}_{\mathbb L^{2,\kappa}_H}\\
&\No{Z^1-Z^2}^2_{\mathbb H^{2,\kappa}_H}+\underset{\mathbb P\in \mathcal P^\kappa_H}{\sup}\mathbb E^\mathbb P\left[\underset{0\leq t\leq T}{\sup}\abs{K_t^{\mathbb P,1}-K_t^{\mathbb P,2}}^2\right]\\
&\leq C\No{\xi^1-\xi^2}_{\mathbb L^{2,\kappa}_H}\left(\No{\xi^1}_{\mathbb L^{2,\kappa}_H}+\No{\xi^1}_{\mathbb L^{2,\kappa}_H}+(\phi^{2,\kappa}_H)^{1/2}+(\psi^{2,\kappa}_H)^{1/2}\right).
\end{align*}
\end{Theorem}

\proof
As in the previous Proposition, we can follow the proof of Lemma $3$ in \cite{ham}, to obtain that there exists a constant $C_\kappa$ depending only on $\kappa$, $T$ and the Lipschitz constant of $\widehat F$, such that for all $\mathbb P$
\begin{equation}
\label{estim2}
\abs{y_t^{\mathbb P,1}-y_t^{\mathbb P,2}}\leq C_\kappa\left(\mathbb E_t^\mathbb P\left[\abs{\xi^1-\xi^2}^\kappa\right]\right)^{\frac1\kappa}.
\end{equation}

Now by definition of our norms, we get from \reff{estim2} and \reff{representationref} that
\begin{equation}
\label{eq.y2}
\No{Y^1-Y^2}_{\mathbb D^{2,\kappa}_H}^2\leq C_\kappa\No{\xi^1-\xi^2}^2_{\mathbb L^{2,\kappa}_H}.
\end{equation}

Applying It\^o's formula to $\abs{Y^1-Y^2}^2$, under each $\mathbb P\in\mathcal P^\kappa_H$, leads to

\begin{align*}
\nonumber\mathbb E^\mathbb P\left[\int_0^T\scriptstyle\abs{\widehat a_t^{\frac12}(Z_t^1-Z_t^2)}^2dt\right]&\leq C\left(\mathbb E^\mathbb P\left[\abs{\xi^1-\xi^2}^2+\int_0^T\abs{Y_t^1-Y_t^2}d\left(K_t^{\mathbb P,1}-K_t^{\mathbb P,2}\right)\right]\right)\\
\nonumber&\hspace{0.9em}+C\mathbb E^\mathbb P\left[\int_0^T\abs{Y_t^1-Y_t^2}\left(\abs{Y_t^1-Y_t^2}+  |\widehat a_t^{\frac12}(Z_t^1-Z_t^2)|\right)dt\right]\\[0.8em]
\nonumber&\leq C\left(\No{\xi^1-\xi^2}_{\mathbb L^{2,\kappa}_H}^2+\No{Y^1-Y^2}^2_{\mathbb D^{2,\kappa}_H}\right)+\frac12\mathbb E^\mathbb P\left[\int_0^T\scriptstyle\abs{\widehat a_t^{\frac12}(Z_t^1-Z_t^2)}^2dt\right]\\
&\hspace{0.9em}+C\No{Y^1-Y^2}_{\mathbb D^{2,\kappa}_H}\left(\mathbb E^\mathbb P\left[\sum_{i=1}^2\left(K_T^i\right)^2\right]\right)^{\frac12}.
\end{align*}

The estimate for $(Z^1-Z^2)$ is now obvious from the above inequality and the estimates of Proposition \ref{estimatesref}. Finally the estimate for the difference of the increasing processes is obvious by definition.
\ep

\section{A direct existence argument}\label{section.3}
We have shown in Theorem \ref{uniqueref} that if a solution exists, it will necessarily verify the representation \reff{representationref}. This gives us a natural candidate for the solution as a supremum of solutions to standard RBSDEs. However, since those BSDEs are all defined on the support of mutually singular probability measures, it seems difficult to define such a supremum, because of the problems raised by the negligible sets. In order to overcome this, Soner, Touzi and Zhang proposed in \cite{stz} a pathwise construction of the solution to a 2BSDE. Let us describe briefly their strategy.

\vspace{0.5em}
The first step is to define pathwise the solution to a standard BSDE. For simplicity, let us consider first a BSDE with a generator equal to $0$. Then, we know that the solution is given by the conditional expectation of the terminal condition. In order to define this solution pathwise, we can use the so-called regular conditional probability distribution (r.p.c.d. for short) of Stroock and Varadhan \cite{str}. In the general case, the idea is similar and consists on defining BSDEs on a shifted canonical space.

\vspace{0.5em}
Finally, we have to prove measurability and regularity of the candidate solution thus obtained, and the decomposition \reff{2bsderef} is obtained through a non-linear Doob-Meyer decomposition. Our aim in this section is to extend this approach to the reflected case.

\subsection{Notations}

For the convenience of the reader, we recall below some of the notations introduced in \cite{stz}.

\vspace{0.5em}
$\bullet$ For $0\leq t\leq T$, denote by $\Omega^t:=\left\{\omega\in C\left([t,T],\mathbb R^d\right),Ê\text{ }w(t)=0\right\}$ the shifted canonical space, $B^t$ the shifted canonical process, $\mathbb P_0^t$ the shifted Wiener measure and $\mathbb F^t$ the filtration generated by $B^t$.

\vspace{0.5em}
$\bullet$ For $0\leq s\leq t\leq T$ and $\omega\in \Omega^s$, define the shifted path $\omega^t\in \Omega^t$
$$\omega^t_r:=\omega_r-\omega_t,\text{ }\forall r\in [t,T].$$

$\bullet$ For $0\leq s\leq t\leq T$ and $\omega\in \Omega^s$, $\widetilde \omega\in\Omega^t$ define the concatenation path $\omega\otimes_t\widetilde \omega\in\Omega^s$ by
$$(\omega\otimes_t\widetilde \omega)(r):=\omega_r1_{[s,t)}(r)+(\omega_t+\widetilde\omega_r)1_{[t,T]}(r),\text{ }\forall r\in[s,T].$$

$\bullet$ For $0\leq s\leq t\leq T$ and a $\mathcal F^s_T$-measurable random variable $\xi$ on $\Omega^s$, for each $\omega \in\Omega^s$, define the shifted $\mathcal F^t_T$-measurable random variable $\xi^{t,\omega}$ on $\Omega^t$ by
$$\xi^{t,\omega}(\widetilde\omega):=\xi(\omega\otimes_t\widetilde \omega),\text{ }\forall \widetilde\omega\in\Omega^t.$$
Similarly, for an $\mathbb F^s$-progressively measurable process $X$ on $[s,T]$ and $(t,\omega)\in[s,T]\times\Omega^s$, the shifted process $\left\{X_r^{t,\omega},r\in[t,T]\right\}$ is $\mathbb F^t$-progressively measurable.

$\bullet$ For a $\mathbb F$-stopping time $\tau$, the r.c.p.d. of $\mathbb P$ (denoted $\mathbb P^\omega_\tau$) is a probability measure on $\mathcal F_T$ such that
$$\mathbb E_\tau^{\mathbb P}[\xi](\omega)=\mathbb E^{\mathbb P^\omega_\tau}[\xi],\text{ for }\mathbb P-a.e.\ \omega.$$

It also induces naturally a probability measure $\mathbb P^{\tau,\omega}$ (that we also call the r.c.p.d. of $\mathbb P$) on $\mathcal F_T^{\tau(\omega)}$ which in particular satisfies that for every bounded and $\mathcal F_T$-measurable random variable $\xi$
$$\mathbb E^{\mathbb P^\omega_\tau}\left[\xi\right]= \mathbb E^{\mathbb P^{\tau,\omega}}\left[\xi^{\tau,\omega}\right].$$

$\bullet$ We define similarly as in Section \ref{section.1} the set $\bar{\mathcal P}^{t}_S$, by restricting to the shifted canonical space $\Omega^t$, and its subset $\mathcal P^{t,\kappa}_H$.

\vspace{0.5em}
$\bullet$ Finally, we define our "shifted" generator
$$\widehat F^{t,\omega}_s(\widetilde\omega,y,z):=F_s(\omega\otimes_t\widetilde\omega,y,z,\widehat a^t_s(\widetilde\omega)), \text{ }\forall (s,\widetilde\omega)\in[t,T]\times\Omega^t.$$

Notice that thanks to Lemma $4.1$ in \cite{stz2}, this generator coincides for $\mathbb P$-a.e. $\omega$ with the shifted generator as defined above, that is to say
$$F_s(\omega\otimes_t\widetilde\omega,y,z,\widehat a_s(\omega\otimes_t\widetilde\omega)).$$

The advantage of the chosen "shifted" generator is that it inherits the uniform continuity in $\omega$ under the $\mathbb L^\infty$ norm of $F$.

\subsection{Existence when $\xi$ is in $\rm{UC_b}(\Omega)$ }

When $\xi$ is in $\rm{UC_b}(\Omega)$, we know that there exists a modulus of continuity function $\rho$ for $\xi$, $F$ and $S$ in $\omega$. Then, for any $0\leq t \leq s \leq T,\ (y,z)\in \left[0,T\right]\times \mathbb R \times \mathbb{R}^d$ and $\omega,\omega'\in \Omega,\ \tilde{\omega}\in\Omega^t$,
\begin{align*}
\left|\xi^{t,\omega}\left(\tilde{\omega}\right)-\xi^{t,\omega'}\left(\tilde{\omega}\right)\right| \leq \rho\left(\No{\omega-\omega'}_t\right) \text{, } \left|\widehat{F}_s^{t,\omega}\left(\tilde{\omega},y,z\right)-\widehat{F}_s^{t,\omega'}\left(\tilde{\omega},y,z\right)\right| \leq \rho\left(\No{\omega-\omega'}_t\right)
\end{align*}
\begin{align*}
\left|S_s^{t,\omega}\left(\tilde{\omega}\right)-S_s^{t,\omega'}\left(\tilde{\omega}\right)\right| \leq \rho\left(\No{\omega-\omega'}_t\right).
\end{align*}

We then define for all $\omega\in\Omega$, $\Lambda\left(\omega\right):=\underset{0\leq s\leq t}{\sup}\Lambda_t\left(\omega\right),$ where
\begin{equation*}
\Lambda_t\left(\omega\right):=\underset{\mathbb P\in\mathcal P^{t}_H}{\sup}\left(\mathbb E^\mathbb P\left[\abs{\xi^{t,\omega}}^2+\int_t^T|\widehat F^{t,\omega}_s(0,0)|^2ds+\left(\underset{t\leq s\leq T}{\sup}(S^{t,\omega}_s)^+\right)^2\right]\right)^{1/2}.
\end{equation*}

Now since $\widehat F^{t,\omega}$ is also uniformly continuous in $\omega$, we have
\begin{equation*}
\Lambda\left(\omega\right)<\infty \text{ for some } \omega\in\Omega \text{ iff it holds for all } \omega\in\Omega.
\end{equation*}

Moreover, when $\Lambda$ is finite, it is uniformly continuous in $\omega$ under the $\mathbb L^{\infty}$-norm and is therefore $\mathcal F_T$-measurable. By Assumption \ref{assump.h2ref}, we have $
\Lambda_t\left(\omega\right)<\infty \text{ for all } \left(t,\omega\right)\in\left[0,T\right]\times\Omega.$

\vspace{0.5em}
To prove existence, we define the following value process $V_t$ pathwise
\begin{equation}
V_t(\omega):=\underset{\mathbb P\in\mathcal P^{t}_H}{\sup}\mathcal Y^{\mathbb P,t,\omega}_t\left(T,\xi\right), \text{ for all } \left(t,\omega\right)\in\left[0,T\right]\times\Omega,
\end{equation}
where, for any $\left(t_1,\omega\right)\in\left[0,T\right]\times\Omega,\ \mathbb P\in\mathcal P^{t_1,\kappa}_H,t_2\in\left[t_1,T\right]$, and any $\mathcal F_{t_2}$-measurable $\eta\in\mathbb L^{2}\left(\mathbb P\right) $, we denote $\mathcal Y^{\mathbb P,t_1,\omega}_{t_1}\left(t_2,\eta\right):= y^{\mathbb P,t_1,\omega}_{t_1}$, where $\left(y^{\mathbb P,t_1,\omega},z^{\mathbb P,t_1,\omega},k^{\mathbb P,t_1,\omega}\right) $ is the solution of the following RBSDE with lower obstacle $S^{t_1,\omega}$ on the shifted space $\Omega^{t_1} $ under $\mathbb P$
\begin{align}\label{eq.bsdeeeeref}
&y^{\mathbb P,t_1,\omega}_{s}=\eta^{t_1,\omega}+\int^{t_2}_{s}\widehat{F}^{t_1,\omega}_{r}\left(y^{\mathbb P,t_1,\omega}_{r},z^{\mathbb P,t_1,\omega}_{r} \right)dr-\int^{t_2}_{s}z^{\mathbb P,t_1,\omega}_{r}dB^{t_1}_{r}+k^{\mathbb P,t_1,\omega}_{t_2}-k^{\mathbb P,t_1,\omega}_{t_1}\\
&\nonumber y_t^{\mathbb P,t_1,\omega}\geq S_t^{t_1,\omega},\text{ }\mathbb P-a.s.\\
&\int_{t_1}^{t_2}\left(y_{s^-}^{\mathbb P,t_1,\omega}-S_{s^-}^{t_1,\omega}\right)dk_s^{\mathbb P,t_1,\omega}=0,\text{ }\mathbb P-a.s.
\label{skoko}
\end{align}

In view of the Blumenthal zero-one law, $\mathcal Y^{\mathbb P,t,\omega}_{t}\left(T,\xi\right) $ is constant for any given $\left(t,\omega\right) $ and $\mathbb P\in\mathcal P^{t,\kappa}_H $. Moreover, since $\omega_0=0 $ for all $\omega\in\Omega $, it is clear that, for the $y^{\mathbb P} $ defined in \reff{rbsde},
\begin{equation*}
\mathcal Y^{\mathbb P,0,\omega}\left(t,\eta\right)=y^{\mathbb P}\left(t,\eta\right) \text{ for all } \omega\in\Omega.
\end{equation*}

\begin{Remark}
We could have defined our candidate solution in another way, using BSDEs instead of RBSDEs, but with a random time horizon. This is based on the link with optimal stopping given by \reff{opt1}. Notice that this approach is similar to the one used by Fabre \cite{fabre} in her PhD thesis when studying 2BSDEs with the $Z$ part of the solution constrained to stay in a convex set. Using this representation as a supremum of BSDEs for a constrained BSDE is particularly efficient, because in general the non-decreasing process added to the solution has no regularity and we cannot obtain stability results. In our case, the two approaches lead to the same result, in particular because the Skorohod condition for the RBSDE allows us to recover stability, as shown in the Lemma below.
\end{Remark}

\begin{Lemma}\label{unifcont}
Let Assumptions \ref{assump.href} and \ref{assump.h2ref} hold and consider some $\xi$ in $\rm{UC_b}(\Omega)$. Then for all $\left(t,\omega\right)\in\left[0,T\right]\times\Omega$ we have $\left|V_t\left(\omega\right)\right|\leq C(1+\Lambda_t\left(\omega\right)) $. Moreover, for all $\left(t,\omega,\omega'\right)\in\left[0,T\right]\times\Omega^2$,
$\left|V_t\left(\omega\right)-V_t\left(\omega'\right)\right|\leq C\rho\left(\No{\omega-\omega'}_t\right) $. Thus, $V_t$ is $\mathcal F_t$-measurable for every $t\in\left[0,T\right]$.
\end{Lemma}

\proof
$\rm{(i)}$ For each $\left(t,\omega\right)\in\left[0,T\right]\times\Omega $ and $\mathbb P\in\mathcal P^{t,\kappa}_H $, let $\alpha$ be some positive constant which will be fixed later and let $\eta\in(0,1)$. By It\^o's formula we have, since $\widehat F$ is uniformly Lipschitz and since by \reff{skoko} $\int_t^Te^{\alpha s}\left(y^{\mathbb P,t,\omega}_{s^-}-S^{t,\omega}_{s^-}\right)dk_s^{\mathbb P,t,\omega}=0$
\begin{align*}
&e^{\alpha t}\abs{y_t^{\mathbb P,t,\omega}}^2+\int_t^Te^{\alpha s}\abs{(\widehat a^t_s)^{1/2}z_s^{\mathbb P,t,\omega}}^2ds\leq e^{\alpha T}\abs{\xi^{t,\omega}}^2+2C\int_t^Te^{\alpha s}\abs{y_s^{\mathbb P,t,\omega}}\abs{\widehat F_s^{t,\omega}(0)}ds\\
&\hspace{0.9em}+2C\int_t^T\abs{y_s^{\mathbb P,t,\omega}}\left(\abs{y_s^{\mathbb P,t,\omega}}+\abs{(\widehat a_s^t)^{1/2}z_s^{\mathbb P,t,\omega}}\right)ds-2\int_t^Te^{\alpha s}y^{\mathbb P,t,\omega}_{s^-}z_s^{\mathbb P,t,\omega}dB^t_s\\
&\hspace{0.9em}+2\int_t^Te^{\alpha s}S_{s^-}^{t,\omega}dk_s^{\mathbb P,t,\omega}-\alpha\int_t^Te^{\alpha s}\abs{y_s^{\mathbb P,t,\omega}}^2ds\\[1em]
&\leq e^{\alpha T}\abs{\xi^{t,\omega}}^2+\int_t^Te^{\alpha s}\abs{\widehat F^{t,\omega}_s(0)}^2ds-2\int_t^Te^{\alpha s}y^{\mathbb P,t,\omega}_{s^-}z_s^{\mathbb P,t,\omega}dB^t_s+\eta\int_t^Te^{\alpha s}\abs{(\widehat a^t_s)^{1/2}z_s^{\mathbb P,n}}^2ds\\[0.3em]
&\hspace{0.9em}+\left(2C+C^2+\frac{C^2}{\eta}-\alpha\right)\int_t^Te^{\alpha s}\abs{y_s^{\mathbb P,t,\omega}}^2ds+2\underset{t\leq s\leq T}{\Sup}e^{\alpha s}(S_s^{t,\omega})^+(k_T^{\mathbb P,t,\omega}-k_t^{\mathbb P,t,\omega}).
\end{align*}

\vspace{0.5em}
Now choose $\alpha$ such that $\nu:=\alpha -2C-C^2-\frac{C^2}{\eta}\geq 0$. We obtain for all $\epsilon>0$
\begin{align}
\nonumber &e^{\alpha t}\abs{y_t^{\mathbb P,t,\omega}}^2+(1-\eta)\int_t^Te^{\alpha s}\abs{(\widehat a_s^t)^{1/2}z_s^{\mathbb P,t,\omega}}^2ds\leq e^{\alpha T}\abs{\xi^{t,\omega}}^2\\
\nonumber&+\int_t^T{e^{\alpha s}{\scriptstyle\abs{\widehat F^{t,\omega}_s(0,0)}^2}}ds+\frac1\epsilon{\scriptstyle\left(\underset{t\leq s\leq T}{\Sup}e^{\alpha s}(S_s^{t,\omega})^+\right)^2}+\epsilon(k_T^{\mathbb P,t,\omega}-k_t^{\mathbb P,t\omega})^2-2\int_t^Te^{\alpha s}y_{s^-}^{\mathbb P,t,\omega}z_s^{\mathbb P,t,\omega}dB_s^t.
\end{align}

Taking expectation yields
\begin{align*}
\abs{y_t^{\mathbb P,t,\omega}}^2+(1-\eta)\mathbb E^\mathbb P\left[\int_t^T\abs{(\widehat a_s^t)^{1/2}z_s^{\mathbb P,t,\omega}}^2ds\right]&\leq C\Lambda_t(\omega)^2
+\epsilon\mathbb E^\mathbb P\left[(k_T^{\mathbb P,t,\omega}-k_t^{\mathbb P,t,\omega})^2\right].
\end{align*}

Now by definition, we also have for some constant $C_0$ independent of $\epsilon$
\begin{align*}
\mathbb E^\mathbb P\left[(k_T^{\mathbb P,t,\omega}-k_t^{\mathbb P,t,\omega})^2\right]&\leq C_0\mathbb E^\mathbb P\left[\abs{\xi^{t,\omega}}^2+\int_t^T\abs{\widehat F^{t,\omega}_s(0,0)}^2ds +\int_t^T\abs{y_s^{\mathbb P,t,\omega}}^2ds\right]\\
&\hspace{0.9em}+\mathbb E^\mathbb P\left[\int_t^T\abs{(\widehat a_s^t)^{1/2}z_s^{\mathbb P,t,\omega}}^2ds\right]\\
&\leq C_0\left(\Lambda_t(\omega)+\mathbb E^\mathbb P\left[\int_t^T\abs{y_s^{\mathbb P,t,\omega}}^2ds+\int_t^T\abs{(\widehat a_s^t)^{1/2}z_s^{\mathbb P,t,\omega}}^2ds\right]\right).
\end{align*}

Choosing $\eta$ small and $\epsilon=\frac{1}{2C_0}$, Gronwall inequality then implies $|y_t^{\mathbb P,t,\omega}|^2\leq C(1+\Lambda_t(\omega)).$ The result then follows by arbitrariness of $\mathbb P$.

\vspace{0.5em}
$\rm{(ii)}$ The proof is exactly the same as above, except that one has to use uniform continuity in $\omega$ of $\xi^{t,\omega}$, $\widehat F^{t,\omega}$ and $S^{t,\omega}$. Indeed, for each $\left(t,\omega\right)\in\left[0,T\right]\times\Omega $ and $\mathbb P\in\mathcal P^{t,\kappa}_H $, let $\alpha$ be some positive constant which will be fixed later and let $\eta\in(0,1)$. By It\^o's formula we have, since $\widehat F$ is uniformly Lipschitz
\begin{align*}
&e^{\alpha t}\abs{y_t^{\mathbb P,t,\omega}-y_t^{\mathbb P,t,\omega'}}^2+\int_t^Te^{\alpha s}\abs{(\widehat a^t_s)^{1/2}(z_s^{\mathbb P,t,\omega}-z_s^{\mathbb P,t,\omega'})}^2ds\leq e^{\alpha T}\abs{\xi^{t,\omega}-\xi^{t,\omega'}}^2\\
&\hspace{0.9em}+2C\int_t^Te^{\alpha s}\abs{y_s^{\mathbb P,t,\omega}-y_s^{\mathbb P,t,\omega'}}\left(\abs{y_s^{\mathbb P,t,\omega}-y_s^{\mathbb P,t,\omega'}}+\abs{(\widehat a_s^t)^{\frac12}(z_s^{\mathbb P,t,\omega}-z_s^{\mathbb P,t,\omega'})}\right)ds\\
&\hspace{0.9em}+2C\int_t^Te^{\alpha s}\abs{y_s^{\mathbb P,t,\omega}-y_s^{\mathbb P,t,\omega'}}\abs{\widehat F^{t,\omega}_s(y_s^{\mathbb P,t,\omega},z_s^{\mathbb P,t,\omega})-\widehat F^{t,\omega'}_s(y_s^{\mathbb P,t,\omega},z_s^{\mathbb P,t,\omega})}ds\\
&\hspace{0.9em}+2\int_t^Te^{\alpha s}(y^{\mathbb P,t,\omega}_{s^-}-y_{s^-}^{\mathbb P,t,\omega'})d(k_s^{\mathbb P,t,\omega}-k_s^{\mathbb P,t,\omega'})-\alpha\int_t^Te^{\alpha s}\abs{y_s^{\mathbb P,t,\omega}-y_s^{\mathbb P,t,\omega'}}^2ds\\
&\hspace{0.9em}-2\int_t^Te^{\alpha s}(y^{\mathbb P,t,\omega}_{s^-}-y_{s^-}^{\mathbb P,t,\omega'})(z_s^{\mathbb P,t,\omega}-z_s^{\mathbb P,t,\omega'})dB^t_s\\
&\leq e^{\alpha T}\abs{\xi^{t,\omega}-\xi^{t,\omega'}}^2+\int_t^Te^{\alpha s}\abs{\widehat F^{t,\omega}_s(y_s^{\mathbb P,t,\omega},z_s^{\mathbb P,t,\omega})-\widehat F^{t,\omega'}_s(y_s^{\mathbb P,t,\omega},z_s^{\mathbb P,t,\omega})}^2ds\\
&\hspace{0.9em}+\left(\scriptstyle2C+C^2+\frac{C^2}{\eta}-\alpha\right)\int_t^Te^{\alpha s}\abs{y_s^{\mathbb P,t,\omega}-y_s^{\mathbb P,t,\omega'}}^2ds+\eta\int_t^Te^{\alpha s}\abs{(\widehat a^t_s)^{\frac12}(z_s^{\mathbb P,t,\omega}-z_s^{\mathbb P,t,\omega'})}^2ds
\\
&\hspace{0.9em}-2\int_t^Te^{\alpha s}(y^{\mathbb P,t,\omega}_{s^-}-y_{s^-}^{\mathbb P,t,\omega'})(z_s^{\mathbb P,t,\omega}-z_s^{\mathbb P,t,\omega'})dB^t_s\\
&\hspace{0.9em}+2\int_t^Te^{\alpha s}(y^{\mathbb P,t,\omega}_{s^-}-y_{s^-}^{\mathbb P,t,\omega'})d(k_s^{\mathbb P,t,\omega}-k_s^{\mathbb P,t,\omega'}).
\end{align*}

By the Skorohod condition \reff{skoko}, we also have
$$\int_t^Te^{\alpha s}(y^{\mathbb P,t,\omega}_{s^-}-y_{s^-}^{\mathbb P,t,\omega'})d(k_s^{\mathbb P,t,\omega}-k_s^{\mathbb P,t,\omega'})\leq \int_t^Te^{\alpha s}(S^{t,\omega}_{s^-}-S_{s^-}^{t,\omega'})d(k_s^{\mathbb P,t,\omega}-k_s^{\mathbb P,t,\omega'}).$$

Now choose $\alpha$ such that $\nu:=\alpha -2C-C^2-\frac{C^2}{\eta}\geq 0$. We obtain for all $\epsilon>0$
\begin{align}
&\nonumber e^{\alpha t}\abs{y_t^{\mathbb P,t,\omega}-y_t^{\mathbb P,t,\omega'}}^2+(1-\eta)\int_t^Te^{\alpha s}\abs{(\widehat a^t_s)^{1/2}(z_s^{\mathbb P,t,\omega}-z_s^{\mathbb P,t,\omega'})}^2ds\\
&\nonumber\leq  e^{\alpha T}\abs{\xi^{t,\omega}-\xi^{t,\omega'}}^2+\int_t^Te^{\alpha s}\abs{\widehat F^{t,\omega}_s(y_s^{\mathbb P,t,\omega},z_s^{\mathbb P,t,\omega})-\widehat F^{t,\omega'}_s(y_s^{\mathbb P,t,\omega},z_s^{\mathbb P,t,\omega})}^2ds\\
\nonumber&\hspace{0.9em}+\frac1\epsilon\left(\underset{t\leq s\leq T}{\Sup}e^{\alpha s}(S_s^{t,\omega}-S_s^{t,\omega'})^+\right)^2+\epsilon(k_T^{\mathbb P,t,\omega}-k_T^{\mathbb P,t,\omega'}-k_t^{\mathbb P,t\omega}+k_t^{\mathbb P,t,\omega'})^2\\
&\hspace{0.9em}-2\int_t^Te^{\alpha s}(y^{\mathbb P,t,\omega}_{s^-}-y_{s^-}^{\mathbb P,t,\omega'})(z_s^{\mathbb P,t,\omega}-z_s^{\mathbb P,t,\omega'})dB^t_s.
\label{bidule3}
\end{align}

The end of the proof is then similar to the previous step, using the uniform continuity in $\omega$ of $\xi$, $F$ and $S$.
\ep

\vspace{0.5em}
Now, we show the same dynamic programming principle as Proposition $4.7$ in \cite{stz2}

\begin{Proposition}\label{progdyn}
Under Assumptions \ref{assump.href}, \ref{assump.h2ref} and for $\xi\in \rm{UC_b}(\Omega)$, we have for all $0\leq t_1<t_2\leq T$ and for all $\omega \in \Omega$
$$V_{t_1}(\omega)=\underset{\mathbb P\in \mathcal P^{t_1,\kappa}_H}{\sup}\mathcal Y_{t_1}^{\mathbb P,t_1,\omega}(t_2,V_{t_2}^{t_1,\omega}).$$
\end{Proposition}

\vspace{0.5em}
The proof is exactly the same as the proof in \cite{stz2}, since we have a comparison Theorem for RBSDEs and since thanks to the paper of Xu and Qian \cite{xuq}, we know that the solution of reflected BSDEs with Lipschitz generator can be constructed via Picard iteration. Given the length of the paper, we omit it. Define now for all $(t,\omega)$, the $\mathbb F^+$-progressively measurable process
\begin{equation}\label{vplus}
V_t^+:=\underset{r\in\mathbb Q\cap(t,T],r\downarrow t}{\overline \lim}V_r.
\end{equation}

We have the following Lemma whose proof is postponed to the Appendix
\begin{Lemma}\label{lemmee}
Under the conditions of the previous Proposition, we have
$$V_t^+=\underset{r\in\mathbb Q\cap(t,T],r\downarrow t}{\lim}V_r,\text{ }\mathcal P_H^\kappa-q.s.$$
and thus $V^+$ is c\`adl\`ag $\mathcal P_H^\kappa-q.s.$.
\end{Lemma}

Proceeding exactly as in Steps $1$ et $2$ of the proof of Theorem $4.5$ in \cite{stz2}, we can then prove that $V^+$ is a strong reflected $\widehat F$-supermartingale. Then, using the Doob-Meyer decomposition proved in the Appendix in Theorem \ref{doobmeyer} for all $\mathbb P$, we know that there exists a unique ($\mathbb P-a.s.$) process $\overline Z^\mathbb P\in\mathbb H^2(\mathbb P)$ and unique non-decreasing c\`adl\`ag square integrable processes $A^\mathbb P$ and $B^\mathbb P$ such that
\begin{itemize}
\item $V_t^+=V_0^+-\int_0^t\widehat F_s(V_s^+,\overline Z_s^\mathbb P)ds+\int_0^t\overline Z^\mathbb P_sdB_s-A_t^\mathbb P-B_t^\mathbb P, \text{ }\mathbb P-a.s., \text{ }\forall \mathbb P\in\mathcal P_H^\kappa.$
\item $V_t^+\geq S_t, \text{ }\mathbb P-a.s. \text{ }\forall \mathbb P\in\mathcal P_H^\kappa.$
\item $\int_0^T\left(V_{t^-}-S_{t^-}\right)dA_t^\mathbb P, \text{ }\mathbb P-a.s., \text{ }\forall \mathbb P\in\mathcal P_H^\kappa.$
\item $A^\mathbb P$ and $B^\mathbb P$ never act at the same time.
\end{itemize}

We then define $K^\mathbb P:=A^\mathbb P+B^\mathbb P$. By Karandikar \cite{kar}, since $V^+$ is a c\`adl\`ag semimartingale, we can define a universal process $\overline Z$ which aggregates the family $\left\{\overline Z^\mathbb P,\mathbb P\in\mathcal P^\kappa_H\right\}$.

\vspace{0.5em}
We next prove the representation \reff{representationref} for $V$ and $V^+$.

\begin{Proposition}\label{prop.repref}
Assume that $\xi\in UC_b(\Omega)$. Under Assumptions \ref{assump.href} and \ref{assump.h2ref}, we have
$$V_t=\underset{\mathbb P^{'}\in\mathcal P_H^\kappa(t,\mathbb P)}{\esup^\mathbb P}\mathcal Y_t^{\mathbb P^{'}}(T,\xi)\text{ and } V_t^+=\underset{\mathbb P^{'}\in\mathcal P_H^\kappa(t^+,\mathbb P)}{\esup^\mathbb P}\mathcal Y_t^{\mathbb P^{'}}(T,\xi), \text{ }\mathbb P-a.s., \text{ }\forall \mathbb P\in\mathcal P_H^\kappa.$$

\end{Proposition}

\proof
The proof for the representations is the same as the proof of proposition $4.10$ in \cite{stz2}, since we also have a stability result for RBSDEs under our assumptions.
\ep

\vspace{0.5em}
Finally, we have to check that the minimum condition \reff{2bsde.minref} holds. Fix $\mathbb P$ in $\mathcal P^\kappa_H$ and $\mathbb P^{'}\in\mathcal P^\kappa_H(t^+,\mathbb P)$. By the Lipschitz property of $F$, we know that there exists bounded processes $\lambda$ and $\eta$ such that
\begin{align}\label{bidulou2}
\nonumber V_t^+-y_t^{\mathbb P^{'}}=&\int_t^T\lambda_s(V_s^+-y_s^{\mathbb P^{'}})ds-\int_t^T\widehat a_s^{1/2}(\overline Z_s-z_s^{\mathbb P^{'}})(\widehat a^{-1/2}_sdB_s-\eta_sds)\\
&+K_T-K_t-k_T^{\mathbb P^{'}}+k_t^{\mathbb P^{'}}.
\end{align}

Then, one can define a probability measure $\mathbb Q^{'}$ equivalent to $\mathbb P^{'}$ such that
$$V_t^+-y_t^{\mathbb P^{'}}=e^{-\int_0^t\lambda_udu}\mathbb E^{\mathbb Q^{'}}_t\left[\int_t^Te^{\int_0^s\lambda_udu}d(K_s-k^{\mathbb P^{'}}_s)\right].$$

Now define the following c\`adl\`ag non-decreasing processes
$$\overline{K}_s:=\int_0^se^{\int_0^u\lambda_rdr}dK_u,\text{ }\overline{k}^{\mathbb P^{'}}_s:=\int_0^se^{\int_0^u\lambda_rdr}dk^{\mathbb P^{'}}_u.$$

By the representation \reff{representationref}, we deduce that the process $\overline{K}-\overline{k}^{\mathbb P^{'}}$
is a $\mathbb Q^{'}$-submartingale. Using Doob-Meyer decomposition and the fact that all the probability measures we consider satisfy the martingale representation property, we deduce as in Step $\rm{(ii)}$ of the proof of Theorem \ref{uniqueref} that this process is actually non-decreasing. Then by definition, this entails that the process $K-k^{\mathbb P^{'}}$ is also non-decreasing.

\vspace{0.5em}
Let us denote $P_t^{\mathbb P^{'}}:=K-k^{\mathbb P^{'}}.$ Returning to \reff{bidulou2} and defining a process $M$ as in Step $\rm{(ii)}$ of the proof of Theorem \ref{uniqueref}, we obtain that
\begin{align*}
V_t^+-y_t^{\mathbb P^{'}}&=\mathbb E^{\mathbb P^{'}}_t\left[\int_t^TM_sdP_s^{\mathbb P^{'}}\right]
\geq \mathbb E^{\mathbb P^{'}}_t\left[\underset{t\leq s\leq T}{\inf}M_s\left(P_T^{\mathbb P^{'}}-P_t^{\mathbb P^{'}}\right)\right].
\end{align*}

Then, we have
\begin{align*}
&\mathbb E^{\mathbb P^{'}}_t\left[P_T^{\mathbb P^{'}}-P_t^{\mathbb P^{'}}\right]=\mathbb E^{\mathbb P^{'}}_t\left[\left(\underset{t\leq s\leq T}{\inf}M_s\right)^{1/3}\left(P_T^{\mathbb P^{'}}-P_t^{\mathbb P^{'}}\right)\left(\underset{t\leq s\leq T}{\inf}M_s\right)^{-1/3}\right]\\
&\leq \left(\mathbb E^{\mathbb P^{'}}_t\left[\underset{t\leq s\leq T}{\inf}M_s\left(P_T^{\mathbb P^{'}}-P_t^{\mathbb P^{'}}\right)\right]\mathbb E^{\mathbb P^{'}}_t\left[\underset{t\leq s\leq T}{\sup}M_s^{-1}\right]\mathbb E^{\mathbb P^{'}}_t\left[\left(P_T^{\mathbb P^{'}}-P_t^{\mathbb P^{'}}\right)^2\right]\right)^{1/3}\\
&\leq C\left(\underset{\mathbb P^{'}\in\mathcal P^{\kappa}_H(t^+,\mathbb P)}{\esup^\mathbb P}\mathbb E^{\mathbb P^{'}}\left[\left(P_T^{\mathbb P^{'}}-P_t^{\mathbb P^{'}}\right)^2\right]\right)^{1/3}\left(V_t^+-y_t^{\mathbb P^{'}}\right)^{1/3}.
\end{align*}

\vspace{0.5em}
Arguing as in Step $\rm{(iii)}$ of the proof of Theorem $\ref{uniqueref}$, we obtain
$$\underset{\mathbb P^{'}\in\mathcal P^{\kappa}_H(t^+,\mathbb P)}{\einf^\mathbb P}\mathbb E^{\mathbb P^{'}}\left[P_T^{\mathbb P^{'}}-P_t^{\mathbb P^{'}}\right]=0,$$
that is to say that the minimum condition \reff{2bsde.minref} is satisfied.

\subsection{Main result}
We are now in position to state the main result of this section

\begin{Theorem}\label{mainref}
Let $\xi\in\mathcal L^{2,\kappa}_H$  and assume that assumptions \ref{assump.href} and \ref{assump.h2ref} hold. Then:\\
1)   There exists a unique solution $(Y,Z)\in\mathbb D^{2,\kappa}_H\times\mathbb H^{2,\kappa}_H$ of the $2\rm{RBSDE}$ \reff{2bsderef}. \\
2) Moreover, if in addition we choose to work under either of the following model of set theory (we refer the reader to \cite{frem} for more details)
\begin{itemize}
\item[\rm{(i)}] Zermelo-Fraenkel set theory with axiom of choice (ZFC) plus the Continuum Hypothesis (CH).
\item[\rm{(ii)}] ZFC plus the negation of CH plus Martin's axiom.
\end{itemize}
Then there exists a unique solution $(Y,Z,K)\in\mathbb D^{2,\kappa}_H\times\mathbb H^{2,\kappa}_H\times\mathbb I^{2,\kappa}_H$ of the $2\rm{RBSDE}$ \reff{2bsderef}.
\end{Theorem}

\proof
The proof of the existence part follows the lines of the proof of Theorem $4.7$ in \cite{stz}, using the estimates of Proposition \ref{estimates2}, so we omit it.  Concerning the fact that we can aggregate the family $\left(K^\mathbb P\right)_{\mathbb P\in\mathcal P^\kappa_H}$, it can be deduced as follows.
\vspace{0.5em}
First, if $\xi\in \rm{UC}_b(\Omega)$, we know, using the same notations as above that our solution verifies
$$V_t^+=V_0^+-\int_0^t\widehat F_s(V_s^+,\overline Z_s)ds+\int_0^t\overline Z_sdB_s-K_t^\mathbb P, \text{ }\mathbb P-a.s., \text{ }\forall \mathbb P\in\mathcal P_H^\kappa.$$

Now, we know from \reff{vplus} that $V^+$ is defined pathwise, and so is the Lebesgue integral $\int_0^t\widehat F_s(V_s^+,\overline Z_s)ds$. In order to give a pathwise definition of the stochastic integral, we would like to use the recent results of Nutz \cite{nutz}. However, the proof in this paper relies on the notion of medial limits, which may or may not exist depending on the model of set theory chosen. They exists in the model (i) above, which is the one considered by Nutz, but we know from \cite{frem} (see statement $22$O(l) page $55$) that they also do in the model (ii). Therefore, provided we work under either one of these models, the stochastic integral $\int_0^t\overline Z_sdB_s$ can also be defined pathwise. We can therefore define pathwise
$$K_t:=V_0^+-V_t^+-\int_0^t\widehat F_s(V_s^+,\overline Z_s)ds+\int_0^t\overline Z_sdB_s,$$
and $K$ is an aggregator for the family $\left(K^\mathbb P\right)_{\mathbb P\in\mathcal P^\kappa_H}$, that is to say that it coincides $\mathbb P-a.s.$ with $K^\mathbb P$, for every $\mathbb P\in\mathcal P^\kappa_H$.

\vspace{0.5em}
In the general case when $\xi\in\mathcal L^{2,\kappa}_H$, the family is still aggregated when we pass to the limit.
\ep

\begin{Remark}
Concerning the models of set theory considered to obtain the aggregation for the family $\left(K^\mathbb P\right)_{\mathbb P\in\mathcal P^\kappa_H}$, even though ZFC is now considered as standard, there are still some controversies about CH. This is the reason why we added the model (ii) above which assumes that CH is false. Consequently, whether one decides to accept this axiom or not, we have a model where the aggregation result holds. Nonetheless, we would like to point out that the Continuum Hypothesis is assumed throughout the books of Dellacherie and Meyer on potential theory (see the last paragraph of page $7$ of \cite{dell}).
\end{Remark}

\section{American Options under volatility uncertainty}\label{section.5}
First let us recall the link between American options and RBSDEs in the classical framework (see \cite{elkarquen} for more details). Let $\mathcal M$ be a standard financial complete market ($d$ risky
asset $S$ and a bond). It is well known that in some constrained
cases the pair wealth-portfolio $(X^\mathbb P,\pi^\mathbb P)$ satisfies: $$
X^\mathbb P_t=\xi+\int_t^Tb(s,X^\mathbb P_s,\pi^\mathbb P_s)ds-\int_t^T\pi^\mathbb P_s\sigma_sdW_s $$ where $W$ is a Brownian motion under the underlying probability measure $\mathbb P$, $b$ is convex
and Lipschitz with respect to $(x,\pi)$. In addition we assume
that the process $(b(t,0,0))_{t\leq T}$ is square-integrable and $(\sigma_t)_{t\leq T}$, the volatility matrix of the
$d$ risky assets, is invertible and its inverse $(\sigma_t)^{-1}$
is bounded. The classical case corresponds to
$b(t,x,\pi)=-r_tx-\pi.\sigma_t\theta_t$, where $\theta_t$ is the
risk premium vector.

\vspace{0.5em}
When the American option is exercised at a
stopping time $\nu\geq t$ the yield is given by $$
\tilde{S}_\nu=S_\nu \textbf{1}_{[\nu<T]}+\xi_T
\textbf{1}_{[\nu=T]}.$$ Let $t$ be fixed and let $\nu \geq t$ be
the exercising time of the contingent claim. Then, since the market is complete, there exists a
unique pair
$(X^\mathbb P_s(\nu,\tilde{S}_\nu),\pi_s^\mathbb P(\nu,\tilde{S}_\nu))=(X^{\mathbb P,\nu}_s,\pi^{\mathbb P,\nu}_s)$
which replicates $\tilde{S}_\nu$, $i.e.$, $$
-dX_s^{\mathbb P,\nu}=b(s,X_s^{\mathbb P,\nu},\pi_s^{\mathbb P,\nu})dt-\pi^{\mathbb P,\nu}_s\sigma_sdW_s,\,\,s\leq \nu;
\,\,X^{\mathbb P,\nu}_\nu=\tilde{S}_\nu.$$ Therefore the price of the
contingent claim is given by $ Y^\mathbb P_t=\underset{\nu\in\mathcal T_{t,T}}{\esup}\text{ }X^\mathbb P_t(\nu,\tilde{S}_\nu).$ Then, the link with RBSDE is given by the following Theorem of \cite{elkarquen}
\begin{Theorem}
There exist $\pi^\mathbb P\in {\mathbb H}^{2}(\mathbb P)$ and a
non-decreasing continuous process $k^\mathbb P$ such that for all $t\in[0,T]$
\begin{align*}
\begin{cases}
&Y^\mathbb P_t=\xi+\int_t^Tb(s,Y^\mathbb P_s,\pi_s^\mathbb P)ds-\int_t^T\pi_s^\mathbb P\sigma_sdW_s+k^\mathbb P_T-k^\mathbb P_t\\
&Y^\mathbb P_t\geq S_t\\
&\int_0^T(Y^\mathbb P_{t}-S_{t})dk_t^\mathbb P=0.
\end{cases}
\end{align*}

Furthermore, the stopping time $D_t^\mathbb P=\inf\{s\geq t,
Y^\mathbb P_s= S_s \}\wedge T$ is optimal after $t$.
\end{Theorem}

Let us now go back to our uncertain volatility framework. The pricing of European contingent claims has already been treated in that context by Avellaneda, L\'evy and Paras in \cite{alp}, Denis and Martini in\cite{denis} with capacity theory and more recently by Vorbrink in \cite{vor} using the G-expectation framework. We still consider a financial market with $d$ risky assets $L^1\ldots L^d$, whose dynamics are given by
$$\frac{dL^i_t}{L^i_t}=\mu^i_tdt+dB^i_t,\text{ }\mathcal P^\kappa_H-q.s.\ \forall i=1\ldots d$$
Then for every $\mathbb P\in \mathcal P^\kappa_H$, the wealth process has the following dynamic
$$X^{\mathbb P}_t=\xi+\int_t^Tb(s,X^{\mathbb P}_s,\pi^{\mathbb P}_s)ds-\int_t^T\pi^{\mathbb P}_sdB_s, \text{ }\mathbb P-a.s.$$

In order to be in our $2$RBSDE framework, we have to assume that $b$ satisfies Assumptions \ref{assump.href} and \ref{assump.h2ref}. In particular, $b$ must satisfy stronger integrability conditions and also has to be uniformly continuous in $\omega$ (when we assume that $\widehat a$ in the expression of $b$ is constant). For instance, in the classical case recalled above, it means that $r$ and $\mu$ must be uniformly continuous in $\omega$, which is the case if for example they are deterministic. We will also assume that $\xi\in\mathcal L^{2,\kappa}_H$. Finally, since $S$ is going to be the obstacle, it has to be uniformly continuous in $\omega$.


\vspace{0.5em}
Following the intuitions in the papers mentioned above, it is natural in our now incomplete market to consider as a superhedging price for our contingent claim
$$Y_t=\underset{\mathbb P^{'}\in\mathcal P^\kappa_H(t^+,\mathbb P)}{\esup^\mathbb P}Y_t^{\mathbb P^{'}},\text{ }\mathbb P-a.s.,\text{ }\forall\mathbb P\in\mathcal P^\kappa_H,$$
where $Y_t^\mathbb P$ is the price at time $t$ of the American contingent claim in the complete market mentioned at the beginning, with underlying probability measure $\mathbb P$. Notice immediately that we do not claim that this price is the superreplicating price in our context, in the sense that it would be the smallest one for which there exists a strategy which superreplicates the American option quasi-surely.

\vspace{0.5em}
The following Theorem is then a simple consequence of the previous one

\begin{Theorem}
There exist $\pi\in {\mathbb H}^{2,\kappa}_H$, a
universal non-decreasing c\`adl\`ag process $K$ such that for all $t\in[0,T]$ and for all $\mathbb P\in\mathcal P^\kappa_H$
\begin{align*}
\begin{cases}
&Y_t=\xi+\int_t^Tb(s,Y_s,\pi_s)ds-\int_t^T\pi_sdB_s+K_T-K_t,\text{ }\mathbb P-a.s.\\
&Y_t\geq S_t, \text{ }\mathbb P-a.s.\\
&K_t-k_t^\mathbb P=\underset{\mathbb P^{'}\in\mathcal P^\kappa_H(t^+,\mathbb P)}{\einf^\mathbb P}\mathbb E^{\mathbb P^{'}}_t\left[K_T-k_T^{\mathbb P^{'}}\right],\text{ }\mathbb P-a.s.
\end{cases}
\end{align*}

Furthermore, for all $\epsilon$, the stopping time $D^\epsilon_t=\inf\{s\geq t,
Y_s\leq S_s+\epsilon,\text{ }\mathcal P^\kappa_H-q.s. \}\wedge T$ is $\epsilon$-optimal after $t$. Besides, for all $\mathbb P$, if we consider the stopping times $D^{\mathbb P,\epsilon}_t=\inf\left\{s\geq t,Y_s^\mathbb P\leq S_s+\epsilon,\text{ }\mathbb P-a.s.\right\}\wedge T$, which are $\epsilon$-optimal for the American options under each $\mathbb P$, then for all $\mathbb P$
\begin{equation}\label{eq.tructruc}
D^\epsilon_t\geq D_t^{\epsilon,\mathbb P},\text{ }\mathbb P-a.s.
\end{equation}
\end{Theorem}

\proof
The existence of the processes is a simple consequence of Theorem \ref{mainref} and the fact that $Y$ is the superhedging price of the contingent claim comes from the representation formula \reff{representationref}. Then, the $\epsilon$-optimality of $D_t^\epsilon$ and the inequality \reff{eq.tructruc} are clear by definition.
\ep

\begin{Remark}
The formula \reff{eq.tructruc} confirms the natural intuition that the smallest optimal time (if exists) to exercise the American option when the volatility is uncertain should be the supremum, in some sense, of all the optimal stopping times for the classical American options for each volatility scenario.
\end{Remark}

\begin{Remark}
As explained in Remark \ref{rem.decomp}, we cannot find a decomposition that would isolate the effects due to the obstacle and the ones due to the second-order. It is not clear neither for the existence of an optimal stopping time.
$D_t=\inf\{s\geq t, Y_{s^-}\leq S_{s^-},\text{ }\mathcal P^\kappa_H-q.s. \}\wedge T$ is not optimal after $t$. Between $t$ and $D_t $, $K^{\mathbb P}$ is reduced to the part related to the second-order. However this part does not verify the minimum condition because it is possible to have $Y_{t^-}>y_{t^-}^\mathbb P=S_{t^-}$, thus the process $k^{\mathbb P}$ is not identically null. For more information on this problem, we would like to refer the reader to the very recent article \cite{etz} which give some specific results for the optimal stopping problem under a non-linear expectation (which roughly corresponds to a 2RBSDE with generator equal to $0$).
\end{Remark}

%
\section*{Acknowledgment}
The authors wish to thank the anonymous referee for all the pertinent remarks she/he made and for having pointed out some references.
\begin{appendix}

\section{Appendix}\label{section.6}

\subsection{Technical proof}
\proof[Proof of Lemma \ref{lemmee}]
For each $\mathbb P$, let $(\bar{\mathcal Y}^{\mathbb P},\bar{\mathcal Z}^\mathbb P)$ be the solution of the BSDE with generator $\widehat F$ and terminal condition $\xi$ at time $T$. We define $\widetilde V^\mathbb P:=V-\bar{\mathcal Y}^\mathbb P.$ Then, $\widetilde V^\mathbb P\geq 0,$ $\mathbb P-a.s.$

\vspace{0.5em}
For any $0\leq t_1< t_2\leq T$, let $(y^{\mathbb P,t_2},z^{\mathbb P,t_2},k^{\mathbb P,t_2}):=(\mathcal Y^\mathbb P(t_2,V_{t_2}),\mathcal Z^\mathbb P(t_2,V_{t_2}),\mathcal K^\mathbb P(t_2,V_{t_2}))$. Since for $\mathbb P-a.e.$ $\omega$, $\mathcal Y^\mathbb P_{t_1}(t_2,V_{t_2})(\omega)=\mathcal Y^{\mathbb P,t_1,\omega}(t_2,V_{t_2}^{t_1,\omega})$, we get from Proposition \ref{progdyn}
$$V_{t_1}\geq y_{t_1}^{\mathbb P,t_2},\text{ }\mathbb P-a.s.$$

\vspace{0.5em}
Denote $\widetilde y_t^{\mathbb P,t_2}:=y_t^{\mathbb P,t_2}-\bar{\mathcal Y}^{\mathbb P}_t,\text{ }\widetilde z_t^{\mathbb P,t_2}:=\widehat a_t^{-1/2}(z_t^{\mathbb P,t_2}-\bar{\mathcal Z}^{\mathbb P}_t).$ Then $\widetilde V_{t_1}^\mathbb P\geq \widetilde y_{t_1}^{\mathbb P,t_2}$ and $(\widetilde y^{\mathbb P,t_2},\widetilde z^{\mathbb P,t_2})$ satisfies the following RBSDE with lower obstacle $S-\bar{\mathcal Y}^{\mathbb P}$ on $[0,t_2]$
$$\widetilde y^{\mathbb P,t_2}_t=\widetilde V_{t_2}^\mathbb P+\int_t^{t_2}f_s^\mathbb P(\widetilde y^{\mathbb P,t_2}_s,\widetilde z^{\mathbb P,t_2}_s)ds-\int_t^{t_2}\widetilde z^{\mathbb P,t_2}_sdW_s^\mathbb P+k^{\mathbb P,t_2}_{t_2}-k^{\mathbb P,t_2}_t,$$
where
$$f_t^\mathbb P(\omega,y,z):=\widehat F_t(\omega,y+\bar{\mathcal Y}^{\mathbb P}_t(\omega),\widehat a_t^{-1/2}(\omega)(z+\bar{\mathcal Z}^{\mathbb P}_t(\omega)))-\widehat F_t(\omega,\bar{\mathcal Y}^{\mathbb P}_t(\omega),\bar{\mathcal Z}^{\mathbb P}_t(\omega)).$$

\vspace{0.5em}
By the definition given in the Appendix, $\widetilde V^\mathbb P$ is a positive weak reflected $f^\mathbb P$-supermartingale under $\mathbb P$. Since $f^\mathbb P(0,0)=0$, we can apply the downcrossing inequality proved in the Appendix in Theorem \ref{down} to obtain classically that for $\mathbb P-a.e.$ $\omega$, the limit
$$\underset{r\in\mathbb Q\cup(t,T],r\downarrow t}{\lim}\widetilde V^\mathbb P_r(\omega)$$
exists for all $t$. Finally, since $\bar{\mathcal Y}^{\mathbb P}$ is continuous, we get the result.
\ep

\vspace{0.5em}
\subsection{Reflected g-expectation}
In this section, we extend some of the results of Peng \cite{pengg} concerning $g$-supersolution of BSDEs to the case of RBSDEs. Let us note that the majority of the following proofs follows straightforwardly from the original proofs of Peng, with some minor modifications due to the added reflection. However, we still provide most of them since, to the best of our knowledge, they do not appear anywhere else in the literature. In the following, we fix a probability measure $\mathbb P$.

\subsubsection{Definitions and first properties}
Let us be given the following objects: a function $g_s(\omega,y,z)$, $\mathbb F$-progressively measurable for fixed $y$ and $z$, uniformly Lipschitz in $(y,z)$, a terminal condition $\xi$ which is $\mathcal F_T$-measurable and in $L^2(\mathbb P)$, and c\`adl\`ag process $V$ and $S$ such that
$$\mathbb E^\mathbb P\left[\int_0^T\abs{g_s(0,0)}^2ds\right]+\mathbb E^\mathbb P\left[\underset{0\leq t\leq T}{\sup}\abs{V_t}^2\right]+\mathbb E^\mathbb P\left[\left(\underset{0\leq t\leq T}{\sup}(S_t)^+\right)^2\right]<+\infty.$$

\vspace{0.5em}
We want to study the following problem. Finding $(y,z,k)\in\mathbb D^2(\mathbb P)\times\mathbb H^2(\mathbb P)\times\mathbb I^2(\mathbb P)$ such that
$$\begin{cases}
\label{rbsde2}
&y_t=\xi+\int_t^Tg_s(y_s,z_s)ds-\int_t^T z_sdW_s+k_{T}-k_{t}+V_T-V_t, \text{ }0\leq t\leq T,\text{ } \mathbb P-a.s.\\
\nonumber&y_t\geq S_t, \text{ }\mathbb P-a.s.\\
\nonumber&\int_0^T\left(y_{s^-}-S_{s^-}\right)dk_s=0, \text{ }\mathbb P-a.s., \text{ }\forall t\in[0,T].
\end{cases}$$

\vspace{0.5em}
We first have a result of existence and uniqueness

\begin{Proposition}
Under the above hypotheses, there exists a unique solution $(y,z,k)\in\mathbb D^2(\mathbb P)\times\mathbb H^2(\mathbb P)\times\mathbb I^2(\mathbb P)$ to the reflected BSDE \reff{rbsde2}.
\end{Proposition}

\vspace{0.5em}
\proof
Consider the following penalized BSDE, whose existence and uniqueness are ensured by the results of Peng \cite{pengg}
$$y^{n}_t=\xi+\int_t^Tg_s(y^{n}_s,z^{n}_s)ds-\int_t^Tz^{n}_sdW_s+k_T^{n}-k_t^{n}+V_T-V_t,$$
where $k_t^{n}:=n\int_0^t(y^{n}_s-S_s)^-ds$.

\vspace{0.5em}
Then, define $\widetilde y_t^{n}:=y_t^{n}+V_t$, $\widetilde \xi:=\xi+V_T$, $\widetilde z_t^{n}:=z_t^{n}$, $\widetilde k_t^{n}:=k_t^{n}$ and $\widetilde g_t(y,z):=g_t(y-V,z)$. We have
$$\widetilde y^{n}_t=\widetilde \xi+\int_t^T\widetilde g_s(\widetilde y^{n}_s,\widetilde z^{n}_s)ds-\int_t^T\widetilde z^{n}_sdW_s+\widetilde k_T^{n}-\widetilde k_t^{n},$$

\vspace{0.5em}
Then, since we know by Lepeltier and Xu \cite{lx}, that the above penalization procedure converges to a solution of the corresponding RBSDE, existence and uniqueness are then simple generalization of the classical results in RBSDE theory.
\ep

\vspace{0.5em}
We also have a comparison theorem in this context

\begin{Proposition}\label{prop.compref}
Let $\xi_1$ and $\xi_2\in L^2(\mathbb P)$, $V^i$, $i=1,2$ be two adapted, c\`adl\`ag processes and $g^i_s(\omega,y,z)$ two functions verifying the above assumptions. Let $(y^i,z^i,k^i)\in \mathbb D^{2}(\mathbb P)\times\mathbb H^{2}(\mathbb P)\times\mathbb I^2(\mathbb P)$, $i=1,2$ be the solutions of the following RBSDEs with lower obstacle $S^i$
$$y_t^i=\xi^i+\int_t^Tg^i_s(y^i_s,z^i_s)ds-\int_t^Tz^i_sdW_s+k_T^i-k_t^i+V_T^i-V_t^i,\text{ }\mathbb P-a.s., \text{ }i=1,2,$$
respectively. If we have $\mathbb P-a.s.$ that $\xi_1\geq \xi_2$, $V^1-V^2$ is non-decreasing, $S^1\geq S^2$, and $g^1_s(y^1_s,z^1_s)\geq g^2_s(y^1_s,z^1_s)$, then it holds $\mathbb P-a.s.$ that for all $t\in [0,T]$
$$Y_t^1\geq Y_t^2.$$

\vspace{0.5em}
Besides, if $S^1=S^2$, then we also have $dK^1\leq dK^2$.
\end{Proposition}

\proof
The first part can be proved exactly as in \cite{elkarkap}, whereas the second one comes from the fact that the penalization procedure converges in this framework.
\ep

 \begin{Remark}
 If we replace the deterministic time $T$ by a bounded stopping time $\tau$, then all the above is still valid.
 \end{Remark}

From now on, we will specialize the discussion to the case where the process $V$ is actually in $\mathbb I^2(\mathbb P)$ and consider the following RBSDE
$$\begin{cases}
\label{rbsde3}
&y_t=\xi+\int_{t\wedge\tau}^\tau g_s(y_s,z_s)ds+V_\tau-V_{t\wedge\tau}+k_\tau-k_{t\wedge\tau}-\int_{t\wedge\tau}^\tau z_sdW_s, \text{ }0\leq t\leq \tau,\text{ } \mathbb P-a.s.\\
\nonumber&y_t\geq S_t, \text{ }\mathbb P-a.s.\\
\nonumber&\int_0^\tau\left(y_{s^-}-S_{s^-}\right)dk_s=0, \text{ }\mathbb P-a.s., \text{ }\forall t\in[0,\tau].
\end{cases}$$

\begin{Definition}
If $y$ is a solution of a RBSDE of the form \reff{rbsde3}, then we call $y$ a reflected $g$-supersolution on $[0,\tau]$. If $V=0$ on $[0,\tau]$, then we call $y$ a reflected $g$-solution.
\end{Definition}

We now face a first difference from the case of non-reflected supersolution. Since in our case we have two increasing processes, if a $g$-supersolution is given, there can exist several increasing processes $V$ and $k$ such that \reff{rbsde3} is satisfied. Indeed, we have the following proposition

\begin{Proposition}\label{uniqueness}
Given $y$ a $g$-supersolution on $[0,\tau]$, there is a unique $z\in\mathbb H^2(\mathbb P)$ and a unique couple $(k,V)\in(\mathbb I^2(\mathbb P))^2$ (in the sense that the sum $k+V$ is unique), such that $(y,z,k,V)$ satisfy \reff{rbsde3}. Besides, there exists a unique quadruple $(y,z,k',V')$ satisfying \reff{rbsde3} such that $k'$ and $V'$ never act at the same time.
\end{Proposition}

\proof
If both $(y,z,k,V)$ and $(y,z^1,k^1,V^1)$ satisfy \reff{rbsde3}, then applying It\^o's formula to $(y_t-y_t)^2$ gives immediately that $z=z^1$ and thus $k+V=k^1+V^1$, $\mathbb P-a.s.$

\vspace{0.5em}
Then, if $(y,z,k,V)$ satisfying \reff{rbsde3} is given, then it is easy to construct $(k',V')$ such that $k'$ only increases when $y_{t^-}=S_{t^-}$, $V'$ only increases when $y_{t^-}>S_{t^-}$ and $V'_t+k'_t=V_t+k_t$, $dt\times d\mathbb P-a.s.$ Moreover, such a couple is unique.
\ep

\begin{Remark}
We give a counter-example to the general uniqueness in the above Proposition. Let $T=2$ and consider the following RBSDE
$$\begin{cases}
\label{rbsde4}
&y_t=-2+2-t+k_2-k_{t}-\int_{t}^2 z_sdW_s, \text{ }0\leq t\leq 2,\text{ } \mathbb P-a.s.\\
\nonumber&y_t\geq -\frac{t^2}{2}, \text{ }\mathbb P-a.s.\\
\nonumber&\int_0^2\left(y_{s^-}+\frac{t^2}{2}\right)dk_s=0, \text{ }\mathbb P-a.s., \text{ }\forall t\in[0,2].
\end{cases}$$

We then have $z=0$, $y_t=1_{0\leq t\leq 1}\left(\frac12-t\right)-\frac{t^2}{2}1_{1<t\leq 2}$ and $k_t=1_{t\geq 1}\frac{t^2-1}{2}$. However, we can also take
$$y'_t=t1_{t\leq 1}+\left(\frac{t^2}{4}+\frac t4+\frac12\right)1_{1<t\leq 2}\text{ and }k'_t=1_{t\geq 1}\left(\frac{t^2}{4}+\frac34t-1\right).$$
\end{Remark}

\vspace{0.5em}
Following Peng \cite{pengg}, this allows us to define

\begin{Definition}
Let $y$ be a supersolution on $[0,\tau]$ and let $(y,z,k,V)$ be the related unique triple in the sense of the RBSDE \reff{rbsde3}, where $k$ and $V$ never act at the same time. Then we call $(z,k,V)$ the decomposition of $y$.
\end{Definition}

\subsubsection{Monotonic limit theorem}
We now study a limit theorem for reflected g-supersolutions, which is very similar to theorems $2.1$ and $2.4$ of \cite{pengg}.

\vspace{0.5em}
We consider a sequence of reflected $g$-supersolutions
$$\begin{cases}
\label{rbsde5}
&y_t^n=\xi^n+\int_{t}^T g_s(y_s^n,z_s^n)ds+V_T^n-V_{t}^n+k_T^n-k_{t}^n-\int_{t}^T z_s^ndW_s, \text{ }0\leq t\leq \tau,\text{ } \mathbb P-a.s.\\
\nonumber&y_t^n\geq S_t, \text{ }\mathbb P-a.s.\\
\nonumber&\int_0^\tau\left(y_{s^-}^n-S_{s^-}\right)dk_s^n=0, \text{ }\mathbb P-a.s., \text{ }\forall t\in[0,T],
\end{cases}$$
where the $V^n$ are in addition supposed to be continuous.

\begin{Theorem}\label{th.monlim}
If we assume that $(y_t^n)$ increasingly converges to $(y_t)$ with
$$\mathbb E^\mathbb P\left[\underset{0\leq t\leq T}{\sup}\abs{y_t}^2\right]<+\infty,$$
and that $(k_t^n)$ decreasingly converges to $(k_t)$, then $y$ is a $g$-supersolution, that is to say that there exists $(z,V)\in\mathbb H^2(\mathbb P)\times \mathbb I^2(\mathbb P)$ such that
$$\begin{cases}
&y_t=\xi+\int_{t}^T g_s(y_s,z_s)ds+V_T-V_{t}+k_T-k_{t}-\int_{t}^T z_sdW_s, \text{ }0\leq t\leq T,\text{ } \mathbb P-a.s.\\
\nonumber&y_t\geq S_t, \text{ }\mathbb P-a.s.\\
\nonumber&\int_0^T\left(y_{s^-}-S_{s^-}\right)dk_s=0, \text{ }\mathbb P-a.s., \text{ }\forall t\in[0,T],
\end{cases}$$

\vspace{0.5em}
Besides, $z$ is the weak (resp. strong) limit of $z^n$ in $\mathbb H^2(\mathbb P)$ (resp. in $\mathbb H^p(\mathbb P)$ for $p<2$) and $V_t$ is the weak limit of $V^n_t$ in $L^2(\mathbb P)$.
\end{Theorem}

\vspace{0.5em}
Before proving the Theorem, we will need the following Lemma
\begin{Lemma}\label{lem}
Under the hypotheses of Theorem \ref{th.monlim}, there exists a constant $C>0$ independent of $n$ such that
$$\mathbb E^\mathbb P\left[\int_0^T\abs{z_s^n}^2ds+(V_T^n)^2+(k_T^n)^2\right]\leq C.$$
\end{Lemma}

\proof
We have
\begin{align}\label{hhh}
\nonumber A_T^n+k_T^n&=y_0^n-y_T^n-\int_0^Tg_s(y_s^n,z_s^n)ds+\int_0^Tz_s^ndW_s\\
&\leq C\left(\underset{0\leq t\leq T}{\sup}\abs{y_t^n}+\int_0^T\abs{z_s^n}ds+\int_0^T\abs{g_s(0,0)}ds +\abs{\int_0^Tz_s^ndW_s}\right).
\end{align}

Besides, we also have for all $n\geq 1$, $y_t^1\leq y_t^n\leq y_t$ and thus $\abs{y_t^n}\leq \abs{y_t^1}+\abs{y_t}$, which in turn implies that
$$\underset{n}{\sup}\mathbb E^\mathbb P\left[\underset{0\leq t\leq T}{\sup}\abs{y_t^n}^2\right]\leq C.$$

Reporting this in \reff{hhh} and using BDG inequality, we obtain
\begin{align}\label{hhhhh}
\nonumber\mathbb E^\mathbb P\left[(V_T^n)^2+(k_T^n)^2\right]&\leq\mathbb E^\mathbb P\left[(V_T^n+k_T^n)^2\right]\\
&\leq C_0\left(1+\mathbb E^\mathbb P\left[\int_0^T\abs{g_s(0,0)}^2ds+\int_0^T\abs{z_s^n}^2ds\right]\right).
\end{align}

Then, using It\^o's formula, we obtain classically for all $\epsilon>0$
\begin{align}\label{hh}
\nonumber\mathbb E^\mathbb P\left[\int_0^T\abs{z_s^n}^2ds\right]&\leq\mathbb E^\mathbb P\left[(y_T^n)^2+2\int_0^Ty_s^ng_s(y_s^n,z_s^n)ds +2\int_0^Ty_{s^-}^nd(V_s^n+k_s^n)\right]\\
&\leq\mathbb E^\mathbb P\left[C\left(1+\underset{0\leq t\leq T}{\sup}\abs{y_t^n}^2\right)+\int_0^T\frac{\abs{z_s^n}^2}{2}ds+\epsilon\left(\abs{V_T^n}^2+\abs{k_T^n}^2\right)\right].
\end{align}

Then, from \reff{hhhhh} and \reff{hh}, we obtain by choosing $\epsilon=\frac{1}{4C_0}$ that
$$\mathbb E^\mathbb P\left[\int_0^T\abs{z_s^n}^2ds\right]\leq C.$$

Reporting this in \reff{hhh} ends the proof.
\ep

\vspace{0.5em}
\proof[Proof of Theorem \ref{th.monlim}]
By Lemma \ref{lem} and its proof we first have
$$\mathbb E^\mathbb P\left[\int_0^T\abs{g_s(y_s^n,z_s^n)}^2ds\right]\leq C\mathbb E^\mathbb P\left[\int_0^T\abs{g_s(0,0)}^2+\abs{y_s^n}^2+\abs{z_s^n}^2ds\right]\leq C.$$

Then we can proceed exactly as in the proof of Theorem $3.1$ in \cite{pengxu}.
\ep

\subsubsection{Doob-Meyer decomposition}
We now introduce the notion of reflected $g$-(super)martingales.

\begin{Definition}
\begin{itemize}
\item[$\rm{(i)}$] A reflected $g$-martingale on $[0,T]$ is a reflected $g$-solution on $[0,T]$.
\item[$\rm{(ii)}$] $(Y_t)$ is a reflected $g$-supermartingale in the strong (resp. weak) sense if for all stopping time $\tau\leq T$ (resp. all $t\leq T$), we have $\mathbb E^\mathbb P[\abs{Y_\tau}^2]<+\infty$ (resp. $\mathbb E^\mathbb P[\abs{Y_t}^2]<+\infty$) and if the reflected $g$-solution $(y_s)$ on $[0,\tau]$ (resp. $[0,t]$) with terminal condition $Y_\tau$ (resp. $Y_t$) verifies $y_\sigma\leq Y_\sigma$ for every stopping time $\sigma\leq\tau$ (resp. $y_s\leq Y_s$ for every $s\leq t$).
\end{itemize}
\end{Definition}

As in the case without reflection, under mild conditions, a reflected $g$-supermartingale in the weak sense corresponds to a reflected $g$-supermartingale in the strong sense. Besides, thanks to the comparison Theorem, it is clear that a $g$-supersolution on $[0,T]$ is also a $g$-supermartingale in the weak and strong sense on $[0,T]$. The following Theorem addresses the converse property, which gives us a non-linear Doob-Meyer decomposition.

\begin{Theorem}\label{doobmeyer}
Let $(Y_t)$ be a right-continuous reflected $g$-supermartingale on $[0,T]$ in the strong sense with
$$\mathbb E^\mathbb P\left[\underset{0\leq t\leq T}{\sup}\abs{Y_t}^2\right]<+\infty.$$

Then $(Y_t)$ is a reflected $g$-supersolution on $[0,T]$, that is to say that there exists a unique triple $(z,k,V)\in\mathbb H^2(\mathbb P)\times\mathbb I^2(\mathbb P)\times\mathbb I^2(\mathbb P)$ such that
$$\begin{cases}
\label{rbsde6}
&Y_t=Y_T+\int_{t}^T g_s(Y_s,z_s)ds+V_T-V_{t}+k_T-k_{t}-\int_{t}^T z_sdW_s, \text{ }0\leq t\leq T,\text{ } \mathbb P-a.s.\\
\nonumber&Y_t\geq S_t, \text{ }\mathbb P-a.s.\\
\nonumber&\int_0^T\left(Y_{s^-}-S_{s^-}\right)dk_s=0, \text{ }\mathbb P-a.s., \text{ }\forall t\in[0,T].\\
\nonumber&\text{$V$ and $k$ never act at the same time.}
\end{cases}$$
\end{Theorem}

We follow again \cite{pengg} and consider the following sequence of RBSDEs
$$\begin{cases}
\label{rbsde7}
&y_t^n=Y_T+\int_{t}^T g_s(y_s^n,z_s^n)ds+n\int_t^T(Y_s-y_s^n)ds+k_T^n-k_{t}^n-\int_{t}^T z_s^ndW_s, \text{ }0\leq t\leq T\\
\nonumber&y_t^n\geq S_t, \text{ }\mathbb P-a.s.\\
\nonumber&\int_0^T\left(y_{s^-}^n-S_{s^-}\right)dk_s^n=0, \text{ }\mathbb P-a.s., \text{ }\forall t\in[0,T],
\end{cases}$$

We have the following Lemma, whose proof is the same a the one of Lemma $3.4$ in \cite{pengg}.
\begin{Lemma}
For all $n$, we have $Y_t\geq y_t^n.$
\end{Lemma}

\vspace{0.5em}
\proof[Proof of Theorem \ref{doobmeyer}]
The uniqueness is due to the uniqueness for reflected $g$-supersolutions proved in Proposition \ref{uniqueness}. For the existence part, we first notice that since $Y_t\geq y_t^n$ for all $n$, by the comparison Theorem for RBSDEs, we have $y_t^n\leq y_t^{n+1}$ and $dk_t^n\geq dk_t^{n+1}$. Therefore they converge monotonically to some processes $y$ and $k$. Besides, $y$ is bounded from above by $Y$. Therefore, all the conditions of Theorem \ref{th.monlim} are satisfied and $y$ is a reflected $g$-supersolution on $[0,T]$ of the form

$$y_t=Y_T+\int_t^Tg_s(y_s,z_s)ds+V_T-V_t+k_T-k_t-\int_t^Tz_sdW_s,$$
where $V_t$ is the weak limit of $V_t^n:=n\int_0^t(Y_s-y_s^n)ds$.

\vspace{0.5em}
From Lemma \ref{lem}, we have
$$\mathbb E^\mathbb P[(V_T^n)^2]=n^2\mathbb E^\mathbb P\left[\int_0^T\abs{Y_s-y_s^n}^2ds\right]\leq C.$$

It then follows that $Y_t=y_t$, which ends the proof.
\ep

\vspace{0.5em}
\subsubsection{Downcrossing inequality}

In this section we prove a downcrossing inequality for reflected $g$-supermartingales in the spirit of the one proved in \cite{cpeng}. We use the same notations as in the classical theory of $g$-martingales (see \cite{cpeng} and \cite{pengg} for instance).

\begin{Theorem}\label{down}
Assume that $g(0,0)=0$. Let $(Y_t)$ be a positive reflected $g$-supermartingale in the weak sense and let $0=t_0<t_1<...<t_i=T$ be a subdivision of $[0,T]$. Let $0\leq<a<b$, then there exists $C>0$ such that $D_a^b[Y,n]$, the number of downcrossings of $[a,b]$ by $\left\{Y_{t_j}\right\}$, verifies
$$\mathcal E^{-\mu}[D_a^b[Y,n]]\leq\frac{C}{b-a}\mathcal E^\mu[Y_0\wedge b],$$
where $\mu$ is the Lipschitz constant of $g$.
\end{Theorem}

\vspace{0.5em}
\proof
Consider
$$\begin{cases}
&y_t^i=Y_{t_i}+\int_{t}^{t_i} ds+\int_t^T(\mu\abs{y_s^i}+\mu\abs{z_s^i}ds+k_T^n-k_{t_j}^n-\int_{t}^{t_i} z_s^idW_s, \text{ }0\leq t\leq t_i,\text{ } \mathbb P-a.s.\\
\nonumber&y_t^i\geq S_t, \text{ }\mathbb P-a.s.\\
\nonumber&\int_0^{t_i}\left(y_{s^-}^i-S_{s^-}\right)dk_s^i=0, \text{ }\mathbb P-a.s., \text{ }\forall t\in[0,t_i].
\end{cases}$$

We define $a_s^i:=-\mu\text{sgn}(z_s^i)1_{t_{j-1}<s\leq t_j}$ and $a_s:=\sum_{i=0}^na_s^i$. Let $\mathbb Q^a$ be the probability measure defined by
$$\frac{d\mathbb Q^a}{d\mathbb P}=\mathcal E\left(\int_0^Ta_sdW_s\right).$$

\vspace{0.5em}
We then have easily that $y_t^i\geq 0$ since $Y_{t_i}\geq 0$ and
$$y_t^i=\underset{\tau\in\mathcal T_{t,t_i}}{\esup}\ \mathbb E_t^{\mathbb Q^a}\left[e^{-\mu(\tau-t)}S_\tau1_{\tau<t_i}+Y_{t_i}e^{-\mu(t_i-t)}1_{\tau=t_i}\right].$$

\vspace{0.5em}
Since $Y$ is reflected $g$-supermartingale (and thus also a reflected $g^{-\mu}$-supermartingale where $g^{-\mu}_s(y,z):=-\mu(\abs{y}+\abs{z})$), we therefore obtain
$$\underset{\tau\in\mathcal T_{t_{i-1},t_i}}{\esup}\mathbb E_{t_{i-1}}^{\mathbb Q^a}\left[e^{-\mu(\tau-t_{i-1})}S_\tau1_{\tau<t_i}+Y_{t_i}e^{-\mu(t_i-t_{i-1})}1_{\tau=t_i}\right]\leq Y_{t_{i-1}}.$$

\vspace{0.5em}
Hence, by choosing $\tau=t_j$ above, we get
$$\mathbb E_{t_{i-1}}^{\mathbb Q^a}\left[Y_{t_i}e^{-\mu(t_i-t_{i-1})}\right]\leq Y_{t_{i-1}},$$
which implies that $(e^{-\mu t_i}Y_{t_i})_{0\leq i\leq n}$ is a $\mathbb Q^a$-supermartingale. Then we can finish the proof exactly as in \cite{cpeng}.
\ep

\end{appendix}

\end{document}